\documentclass[12pt,reqno]{article}
\usepackage{amsmath,amsthm,amsfonts,amssymb,amscd,amstext}
\usepackage[ansinew]{inputenc}
\usepackage[dvips]{graphicx}
\usepackage{psfrag,euscript}
\usepackage{a4wide}

\usepackage{pxfonts}   
\numberwithin{equation}{section}

\theoremstyle{plain}
\newtheorem{maintheorem}{Theorem}
\newtheorem{maincorollary}[maintheorem]{Corollary}

\newtheorem{theorem}{Theorem}[section]
\newtheorem{proposition}[theorem]{Proposition}
\newtheorem{corollary}[theorem]{Corollary}
\newtheorem{lemma}[theorem]{Lemma}

\theoremstyle{definition}

\newcommand{\RR}{{\mathbb R}}

\newcommand{\NN}{{\mathbb N}}
\newcommand{\ZZ}{{\mathbb Z}}
\newcommand{\FF}{{\mathbb F}}
\newcommand{\EE}{{\mathbb E}}

\newcommand{\sS}{{\mathbb S}}

\newcommand{\vfi}{\varphi}

\newcommand{\de}{\delta}

\newcommand{\diam}{\operatorname{diam}}
\renewcommand{\epsilon}{\varepsilon}
\newcommand{\dist}{\operatorname{dist}}

\newcommand{\Leb}{\operatorname{Leb}}

\newcommand{\co}{\operatorname{co}}

\newcommand{\supp}{\operatorname{supp}}

\newcommand{\cE}{\EuScript{E}}

\newcommand{\D}{\EuScript{D}}
\newcommand{\cP}{\EuScript{P}}

\newcommand{\U}{\EuScript{U}}

\newcommand{\cC}{\EuScript{C}}
\newcommand{\cS}{\EuScript{S}}
\newcommand{\M}{\EuScript{M}}

\title
{Large deviations bound for semiflows over a non-uniformly expanding base}

\author{Vítor Araújo (\texttt{vitor.araujo@im.ufrj.br \text{and} vdaraujo@fc.up.pt}) 
\thanks{The author was partially supported by
    CNPq-Brazil and FCT-Portugal through CMUP and
    POCI/MAT/61237/2004. }
\\
Instituto de Matemática,
Universidade Federal do Rio de Janeiro
\\
C. P. 68.530,
  21.945-970 Rio de Janeiro, RJ-Brazil 
\\
\emph{and} 
\\
Centro de Matemática da Universidade do Porto
\\
Rua do Campo Alegre 687, 
4169-007 Porto, Portugal.}

\date{\today}

\begin{document}

\maketitle


\bigskip

\noindent
\textbf{\textup{2000} Mathematics Subject Classification:}
Primary: 37D25; Secondary:
37D35, 37D50, 37C40.

\begin{abstract}
  We obtain a exponential large deviation upper bound for
  continuous observables on suspension semiflows over a
  non-uniformly expanding base transformation with non-flat
  singularities or criticalities, where the roof function
  defining the suspension behaves like the logarithm of the
  distance to the singular/critical set of the base map.
  That is, given a continuous function we consider its space
  average with respect to a physical measure and compare
  this with the time averages along orbits of the semiflow,
  showing that the Lebesgue measure of the set of points
  whose time averages stay away from the space average tends
  to zero exponentially fast as time goes to infinity.

  The arguments need the base transformation to exhibit
  exponential slow recurrence to the singular set which, in
  all known examples, implies exponential decay of
  correlations.

  Suspension semiflows model the dynamics of flows admitting
  cross-sections, where the dynamics of the base is given by
  the Poincar\'e return map and the roof function is the
  return time to the cross-section. The results are
  applicable in particular to semiflows modeling the
  geometric Lorenz attractors and the Lorenz flow, as well
  as other semiflows with multidimensional non-uniformly
  expanding base with non-flat singularities and/or
  criticalities under slow recurrence rate conditions to
  this singular/critical set.  We are also able to obtain
  exponentially fast escape rates from subsets without full
  measure.
\end{abstract}


\section{Introduction}
\label{sec:introd}

The statistical viewpoint on Dynamical Systems provides some
of the main tools available for the global study of the
asymptotic behavior of transformations or flows. One of
the main concepts introduced is the notion of
\emph{physical} (or \emph{Sinai-Ruelle-Bowen}) measure for a
flow or a transformation. An invariant probability measure
$\mu$ for a flow $X^t$ on a compact manifold is a physical
probability measure if the points $z$ satisfying for all
continuous functions $\psi$
\begin{align*}
  \lim_{t\to+\infty}\frac1t\int_0^t \psi\big(X^s(z)\big) \,
  ds = \int\psi\,d\mu,
\end{align*}
form a subset with positive volume (or positive
Lebesgue measure) on the ambient space. These time averages
are in principle physically observable if the flow models
a real world phenomenon admitting some measurable features.

For systems admitting such invariant probability measures it
is natural to consider the rate of convergence of the time
averages to the space average, given by the volume of the
subset of points whose time averages stay away from the
space average by a prescribed amount up to some evolution
time.  This rate is closely related to the so-called
thermodynamical formalism first developed for (uniformly)
hyperbolic diffeomorphisms, borrowed from statistical
mechanics by Bowen, Ruelle and Sinai (among others, see e.g.
\cite{Bo75,BR75,Ru89,ruelle2004,ellis06,BDV2004}). These
authors systematically studied the construction and properties
of physical measures for (uniformly) hyperbolic
diffeomorphisms and flows. Such measures for non-uniformly
hyperbolic maps and flows where obtained more recently
\cite{PS82,CT88,BeY92,BeY93,Al00}.

The probabilistic properties of physical measures are an
object of intense study, see e.g.
\cite{BR75,kifer90,Yo90,Yo98,BeY99,AA03,alves-araujo2004,
  alves-luzzatto-pinheiro2005,
  gouezel,arbieto-matheus2006,araujo-pacifico2006}.  The
main insight behind these efforts is that the family $\{
\psi\circ X^t\}_{t>0}$ should behave asymptotically in many
respects just like a i.i.d. random variable.

The study of suspension (or special) flows is motivated by
modeling a flow admitting a cross-section. Such flow is
equivalent to a suspension semiflow over the Poincar\'e
return map to the cross-section with roof function given by
the return time function for the points in the
cross-section. This is one of the main technical tools in
the ergodic theory of Axiom A (or uniformly hyperbolic)
flows developed by Bowen and Ruelle~\cite{BR75}, enabling
them to pass from this type of flow to a suspension flow
over a shift transformation with finitely many symbols and
bounded roof function. Then the properties of the base
transformation are used to deduce many results for the
suspension flow, which are then passed to the original flow.

Recently, based on the breakthrough of Dolgopyat
\cite{Do98}, this kind of modeling provided results on the
rate of decay of correlations for certain
flows~\cite{AvGoYoc} based on the rate of decay of
correlations for suspension semiflows \cite{BaVal2005}.
General results on the existence and some statistical
properties of physical measures for singular-hyperbolic
attractors for three-dimensional flows \cite{APPV} as well
as their sensitive dependence on initial conditions were
also obtained using this standard technique.  Moreover the
classical Lorenz flow \cite{Lo63} was shown to be equivalent
to a geometric Lorenz flow by Tucker \cite{Tu99} and so it
can be modeled by a suspension semiflow over a non-uniformly
hyperbolic transformation with unbounded roof function.
Using these ideas it was recently obtained \cite{LMP05} that
the physical measure for the Lorenz attractor is mixing.

Here we extend part of the results on large deviation rates
of Kifer \cite{kifer90} (see also Waddington \cite{wadd96})
from the uniformly hyperbolic setting to semiflows over
non-uniformly expanding base dynamics and unbounded roof
function. These special flows model non-hyperbolic flows,
like the Lorenz flow, exhibiting equilibria accumulated by
regular orbits. We use the properties of non-uniformly
expanding transformations, especially the large deviation
bound obtained in \cite{araujo-pacifico2006}, to deduce a
large deviation bound for the suspension semiflow reducing
the estimate of the volume of the deviation set to the
volume of a certain deviation set for the base
transformation.  More precisely, if we set $\epsilon>0$ as
an error margin and consider
\[
B_t=\Big\{
z: \Big|
\frac1t\int_{0}^{t}
  \psi\big( X^t(z) \big) 
-
\int\psi\,d\mu
\Big|>\epsilon
\Big\}
\]
then we are able to provide conditions under which the
Lebesgue measure of $B_t$ decays to zero exponentially fast,
i.e. weather there are constants $C,\xi>0$ such that
\begin{align*}
\Leb\big( B_t\big) \le C e^{-\xi t}
\quad\mbox{for all}\quad t>0.
\end{align*}
The values of $C,\xi>0$ above depend on $\epsilon,\psi$ and
on global invariants for the base transformation $f$, such
as the metric entropy and the pressure function of $f$ with
respect to the physical measures of $f$ and a certain
observable constructed from $\psi$ and $X^t$, as detailed in
the next section. Having this it is not difficult to deduce
exponential escape rates from subsets of the semiflow.

In order to be able to apply this bound to Lorenz flows, it
is necessary to allow the roof function of the suspension
flows to be unbounded near the singularities of the base
dynamical system. This in turn imposes some restrictions on
the admissible base dynamics, expressed as a slow recurrence
rate to the singular set and uniqueness of equilibrium
states with respect to the logarithm of the Jacobian of the
map. However no cohomology condition on the roof function
are needed, while this is essential to obtain fast decay of
correlations in \cite{dolgopyat98,Mel06,FMT}.

We present several semiflows with non-uniformly expanding
base transformations satisfying all our
conditions, including one-dimensional
piecewise expanding maps with \emph{Lorenz-like}
singularities and quadratic maps but also multidimensional
examples. This demanded the detailed study of recurrence rates
to the singular set, the study of large deviation bounds for
unbounded observables over non-uniformly expanding
transformations, and an entropy formula for non-uniformly
expanding maps with singularities (which might be of
independent interest). Now we give the precise statement of
the results.

\subsection{Statement of the results}
\label{sec:statem-results}

Denote by $\|\cdot\|$ a Riemannian norm on the compact
boundaryless manifold $M$, by $\dist$ the induced distance
and by $\Leb$ the corresponding Riemannian volume form,
which we call \emph{Lebesgue measure} or \emph{volume}. We
assume $\Leb$ to be normalized: $\Leb(M)=1$.

Given a $C^2$ local diffeomorphism (H\"older-$C^1$ is
enough, see below) $f:M\setminus\cS\to M$ outside a volume
zero non-flat singular set, let $X^t:M_r\to M_r$ be a
semiflow with roof function $r:M\setminus\cS\to\RR$ over the
base transformation $f$, as follows.  Set $M_r=\{ (x,y)\in
M\times[0,+\infty): 0\le y < r(x) \}$.  For $x=x_0\in M$
denote by $x_n$ the $n$th iterate $f^n(x_0)$ for $n\ge0$.
Denote $S_n \vfi(x_0) = S_n^f \vfi(x_0) = \sum_{j=0}^{n-1}
\vfi( x_j )$ for $n\ge1$ and for any given real function
$\vfi$ in what follows. Then for each pair $(x_0,s_0)\in
X^r$ and $t>0$ there exists a unique $n\ge1$ such that $S_n
r(x_0) \le s_0+ t < S_{n+1} r(x_0)$ and we define
\begin{align*}
   X^t(x_0,s_0) = \big(x_n,s_0+t-S_n r(x_0)\big).
\end{align*}

The non-flatness of the singular set $\cS$ is an extension
to arbitrary dimensions of the notion of non-flat singular
set from one-dimensional dynamics \cite{MS93} and means that
$f$ \emph{behaves like a power of the distance to the
  singular set}. More precisely there are constants $B>1$
and $0<\beta<1$ for which
 \begin{itemize}
 \item[(S1)]
\hspace{.1cm}$\displaystyle{\frac{1}{B} \dist(x,\cS)^{\beta}\leq
\frac{\|Df(x)v\|}{\|v\|}\leq B \dist(x,\cS)^{-\beta}}$;
 \item[(S2)]
\hspace{.1cm}$\displaystyle{\left|\log\|Df(x)^{-1}\|-
\log\|Df(y)^{-1}\|\:\right|\leq
B\frac{ \dist(x,y)}{\dist(x,\cS)^{\beta}}}$;
 \item[(S3)]
\hspace{.1cm}$\displaystyle{\left|\log|\det Df(x)^{-1}|-
\log|\det Df(y)^{-1}|\:\right|\leq
B\frac{\dist(x,y)}{\dist(x,\cS)^{\beta}}}$;
 \end{itemize}
 for every $x,y\in M\setminus \cS$ with
 $\dist(x,y)<\dist(x,\cS)/2$ and $v\in T_x M\setminus\{0\}$.
We also assume an extra condition related to the geometry of
$\cS$. This ensures that the Lebesgue measure of neighborhoods
$\cS$ is comparable to a power of the
distance to $\cS$, that is there exists $C_\kappa,\kappa>0$
such that for all small $\rho>0$
\begin{itemize}
\item[(S4)]
$\Leb\{x\in M: \dist(x,\cS)<\rho\}\le
C_\kappa\cdot\rho^\kappa.$
\end{itemize}
The singular set $\cS$ contains those points $x$ where $f$
is either not defined, is discontinuous, not differentiable
or else $Df(x)$ is non-invertible (that is $\cS$ contains
the set $\cC$ of critical points of $f$).
Note that condition (S4) is satisfied in the particular case
when $\cS$ is a compact submanifold of $M$, where
$\kappa=\dim(M)-\dim(\cS)$. It is also satisfied for
$M=\sS^1$ and $\cS$ is a denumerable infinite subset with
finitely many accumulation points, with $\kappa=1$. In particular
this holds for a piecewise expanding map over the interval
or the circle with finitely many domains of monotonicity.

We say that $f$ is \emph{non-uniformly expanding} if there
exists $c>0$ such that
\begin{align}
  \label{eq:NUE}
  \limsup_{n\to+\infty}\frac1n S_n\psi(x) \le -c
\quad\mbox{where}\quad
\psi(x)=\log\big\| Df(x)^{-1} \big\|,
\end{align}
for Lebesgue almost every $x\in M$. This condition implies
in particular that all the lower Lyapunov exponents of the
map $f$ are strictly positive Lebesgue almost everywhere.

Let $\Delta_\delta(x)=\big| \log d_{\delta}(x,\cS) \big|$ be
the \emph{smooth $\delta$-truncated logarithmic distance} from $x\in M$
to $\cS$,
i.e. $\Delta_\delta(x)$ is non-negative and continuous away
from $\cS$, identically zero $2\delta$-away from $\cS$, and
equal to $-\log\dist(x,\cS)$ when $\dist(x,\cS)\le\delta$.

We say that $f$ has \emph{exponentially slow recurrence to
  the singular set $\cS$} if for every $\epsilon>0$ there
exists $\delta>0$ such that
\begin{align}\label{eq:expslowrecurrence}
  \limsup_{n\to+\infty}\frac1n\log\Leb\left\{x\in M:
    \frac1n S_n\Delta_\delta(x)>\epsilon\right\}<0.
\end{align}
Condition~\eqref{eq:expslowrecurrence} implies that
$S_n\Delta_\delta/n\to0$ in measure, i.e.  for every
$\epsilon>0$ there exists $\delta>0$ such that
\begin{align}
  \label{eq:SlowApprox}
  \limsup_{n\to\infty} \frac1n S_n\Delta_\delta(x)\le \epsilon
\end{align}
for Lebesgue almost every $x\in M$. We say that a map $f$
satisfying \eqref{eq:SlowApprox} has \emph{slow recurrence
  to $\cS$}.

These notions were presented in~\cite{ABV00} and in
\cite{ABV00,Ze03} the following result on existence of
finitely many absolutely continuous measures was obtained.
\begin{theorem}
  \label{thm:abv}
  Let $f:M\to M$ be a $C^2$ local diffeomorphism outside a
  singular set $\cS$. Assume that $f$ is non-uniformly
  expanding with slow recurrence to $\cS$.  Then there are
  finitely many ergodic absolutely continuous (in particular
  \emph{physical }or \emph{Sinai-Ruelle-Bowen})
  $f$-invariant probability measures $\mu_1,\dots,\mu_k$
  whose basins cover the manifold Lebesgue almost
  everywhere, that is $ B(\mu_1)\cup\dots\cup B(\mu_k) =
  M,\quad \Leb-\bmod0$.  Moreover the support of each
  measure contains an open disk in $M$.
\end{theorem}
Here the \emph{basin} of an invariant probability measure
$\mu$ is the subset of points $x\in M$ such that
$\lim_{n\to\infty}\frac1n\sum_{j=0}^{n-1} \delta_{f^j(x)} =
  \mu$ in the weak$^*$ topology.

  Large deviation bounds are usually related to measure
  theoretic entropy and to equilibrium states.  We denote by
  $\M_f$ the family of all invariant probability measures
  with respect to $f$. Let $J=|\det Df|$. We say that
  $\mu\in\M_f$ is an \emph{equilibrium state} with respect
  to the potential $\log J$ if $h_\mu(f)=\mu(\log J)$, that
  is if \emph{$\mu$ satisfies the Entropy Formula}. We
  denote by $\EE$ the subset of $\M_f$ consisting of all
  equilibrium states for $f$. It is not difficult to see
  (Section~\ref{sec:isolat-physic-measur} for more details)
  that each physical measure provided by
  Theorem~\ref{thm:abv} belongs to $\EE$.

Another standing assumption on $f$ is that \emph{the set
  $\EE$ is formed by a unique $f$-invariant absolutely
  continuous probability measure} (see
Section~\ref{sec:exampl-applic} for sufficient conditions
for this to occur and for examples of application).

We denote by $\nu=\mu\ltimes\Leb^1$ the natural $X^t$-invariant
extension of $\mu$ to $M_r$ and by $\lambda$ the natural
extension of $\Leb$ to $M_r$,
i.e. $\lambda=\Leb\ltimes\Leb^1$, where $\Leb^1$ is
one-dimensional Lebesgue measure on $\RR$: for any subset
$A\subset M_r$
\begin{align*}
  \nu(A)=\frac1{\mu(r)}\int d\mu(x)
\int_0^{r(x)}\!\!\!\! ds \,\chi_A(x,s) 
\quad\text{and}\quad
  \lambda(A)=\frac1{\Leb(r)}\int d\Leb(x)
\int_0^{r(x)}\!\!\!\! ds \,\chi_A(x,s). 
\end{align*}

We say that a function $\vfi:M\setminus\cS\to\RR$ has
\emph{logarithmic growth near $\cS$} if there exists
$K=K(\vfi)>0$ such that 
\begin{align}\label{eq:log-growth}
  |\vfi|\chi_{B(\cS,\delta)}\le K\cdot\Delta_\delta 
  \text{   for all small enough   }
  \delta>0.
\end{align}
We also say that $f$ is a \emph{regular map} if for
$E\subset M$ such that $\Leb(E)=0$, then
$\Leb\big(f^{-1}(E)\big)=0$.

\begin{maintheorem}
  \label{mthm:LDNUEflow}
  Let $X^t$ be a suspension semiflow over a non-uniformly
  expanding transformation $f$ on the base $M$ which
  exhibits exponentially slow recurrence to the 
  singular set, where the roof function
  $r:M\setminus\cS\to\RR$ has logarithmic growth near $\cS$.
  Assume that $f$ is a regular map and that the set $\EE$ of
  equilibrium states is formed by a single measure $\mu$.
  Let $\psi:M_r\to\RR$ be a continuous function.  Then 
\begin{align}
  \label{eq:flowDeviation}
\limsup_{T\to\infty}\frac1T\log\lambda\left\{
z\in M_r: \left|
\frac1T\int_0^T \psi\left(X^t(z)\right)\, dt - \nu(\psi)
\right|>\epsilon 
\right\}<0.
\end{align}
\end{maintheorem}

\subsection{Escape rates}
\label{sec:escape-rates}

Let $K\subset M_r$ be a compact subset.  Given $\epsilon>0$
we can find an open set $W\supset K$ contained in $M_r$ and
a continuous bump function $\vfi:M_r\to\RR$ such that
$\Leb(W\setminus K) < \epsilon$ with $0\le\vfi\le1$,
$\vfi\mid K\equiv 1$ and $\vfi\mid (M\setminus W) \equiv
0$. Then we get for $n\ge1$
\begin{align}
  \label{eq:escape}
\left\{
x\in K : X^t(x)\in K, 0<t<T
\right\}
\subset
\left\{
x\in M : \frac1T\int_0^T\vfi\big(X^t(x)\big)dt \ge1 
\right\}
\end{align}
and so we deduce the following 
using the estimate from Theorem~\ref{mthm:LDNUEflow}.
\begin{maincorollary}
  \label{mcor:escaperate}
  Let $X^t$ be a suspension semiflow over a non-uniformly
  expanding transformation $f$ on the base $M$ in the same
  setting as in Theorem~\ref{mthm:LDNUEflow}.  Let $K$ be a
  compact subset of $M_r$ such that $\nu(K)<1$.  Then 
\[
\limsup_{T\to+\infty}
\frac1T \log \lambda\Big( \left\{ x\in K :
  X^t(x)\in K, 0<t<T \right\} \Big) < 0.
\]
\end{maincorollary}

\subsection{Lorentz and Geometric Lorenz flows}
\label{sec:lorentz-geometr-lore}

\begin{figure}[htb]
  \centering
  \psfrag{e}{$\ell$}\psfrag{y}{$f(x)$}\psfrag{S}{$S$}
  \includegraphics[width=10cm]{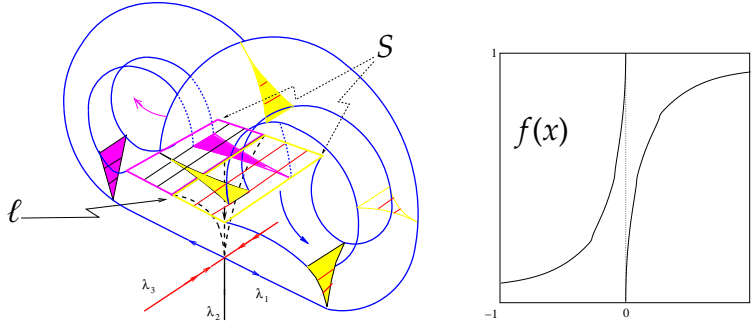}
  \caption{The geometric Lorenz flow and the associated
    one-dimensional piecewise expanding map}
  \label{fig:geometr-lorenz-flow}
\end{figure}

The Lorenz equations
\begin{align}\label{eq:lorentz}
  \dot x = 10(y-x), \quad
  \dot y = 28x -y -xz, \quad
  \dot z = xy - 8z/3
\end{align}
were presented by Lorenz \cite{Lo63} in 1963 as a simplified
model of convection of the Earth's atmosphere. It turned out
that these equations became one of the main models showing
the presence of chaotic dynamics in apparently simple
systems. More recently Tucker \cite{Tu99,Tu2} with a
computer assisted proof showed that equations
\eqref{eq:lorentz} and similar equations with nearby
parameters define a geometric Lorenz flow, i.e. a
three-dimensional flow $X^t$ in $\RR^3$ with a hyperbolic
singularity at the origin admitting a neighborhood $U$ (a
\emph{trapping region}) such that $\overline{X^t(U)}\subset
U$ for all $t>0$ satisfying:
\begin{enumerate}
\item the attracting set $\Lambda=\cap_{t>0}X^t(U)$ contains
  the singularity at $0$;
\item $\Lambda$ contains a dense orbit;
\item there exists a square $S=[-1,1]\times[-1,1]\times\{1\}$
  which is a cross-section for $\Lambda\setminus\{0\}$, that
  is for every $w\in\Lambda\setminus\{0\}$ there exists
  $t>0$ such that $X^t(w)\in S$;
\item the Poincar\'e first return map to $S$ given by
  $R:S\setminus\ell\to S$ is $C^2$ and contracts distances
  exponentially on the $y$ direction, where
  $\ell=\{0\}\times[-1,1]\times\{1\}$ is the singular line,
  so each segment $S\cap\{x=\mathrm{const}\}$ is contained
  in a stable manifold. Moreover in general this
  one-dimensional and co-dimension one foliation of the
  cross-section $S$ defines a projection $P$ along leaves
  which is $C^{1+\alpha}$ for some $\alpha>0$;
\item the one-dimensional map
  $f:[-1,1]\setminus\{0\}\to[-1,1]$ obtained from $R$
  quotienting out the stable manifolds is a piecewise
  expanding map with singularities known as
  \emph{Lorenz-like map}, which is in the setting of the
  class of examples detailed in
  Subsection~\ref{sec:semifl-over-non-1};
\item the roof function $\tau(w)$ for $w\in S$
  is Lebesgue integrable over $S$ and has logarithmic
  growth near the singular line $\ell$.
\end{enumerate}
It is well known that the attractor of the geometric Lorenz
flows (and the attractor for the Lorenz equations after the
results of Tucker already mentioned) supports a unique
ergodic physical measure $\mu$ (for more details on this
construction see e.g. \cite{Vi97b}).
Figure~\ref{fig:geometr-lorenz-flow} gives a visual idea of
the geometric Lorenz flow. The reader should consult
\cite{Gu76,GW79,robinson2004} for proofs of the properties
stated above and more details on the construction of such
flows.
Using $\tau$ as a roof function over the base dynamics given
by $R$ we see that the dynamics of a geometric Lorenz flow on
$U$ is equivalent to a suspension semiflow over $R$ with
roof function $\tau$. In addition the uniform contraction
along the leaves of the foliation $\{y=\mathrm{const}\}$
together with the uniform expansion of the one-dimensional
map $f$ enables us to use Theorem~\ref{mthm:LDNUEflow} to
deduce

\begin{maincorollary}
  \label{mcor:LDLorenz}
  Let $X^t$ be a flow on $\RR^3$ exhibiting a Lorenz or a
  geometric Lorenz attractor with trapping region
  $U$. Denoting by $\Leb$ the normalized restriction of the
  Lebesgue volume measure to $U$, $\psi:U\to\RR$ a
  bounded continuous function and $\mu$ the unique physical
  measure for the attractor, then for any given $\epsilon>0$
  \begin{align*}
      \limsup_{T\to\infty}\frac1T\log\Leb\left\{
        z\in U: \left|
          \frac1T\int_0^T \psi\left(X^t(z)\right)\, dt - \mu(\psi)
        \right|>\epsilon 
      \right\}<0,
  \end{align*}
and consequently for any compact $K\subset U$ such that
$\mu(K)<1$ we also have
\begin{align*}
\limsup_{T\to+\infty}\frac1T \log \Leb\Big( \left\{ x\in K :
  X^t(x)\in K, 0<t<T \right\} \Big) < 0.
\end{align*}
\end{maincorollary}

\subsection{Comments and organization of the paper}
\label{sec:organiz-paper}

We note that the smoothness assumption needed for our
arguments is only $C^{1+\alpha}$ for some
$\alpha\in(0,1)$. Therefore the $C^2$ condition on $f$ in
the statements of results can be relaxed to $C^{1+\alpha}$
throughout.

Kifer~\cite{kifer90} together with Newhouse~\cite{KiNw91}
obtain sharp large deviations bounds both from above and
from below for uniformly partially hyperbolic attractors for
flows and for Axiom A flows, through an estimate of the
volume growth of images of balls under the action of the
flow near the attractor (``volume lemma'', see also
\cite{BR75} and \cite{Bo75}). Moreover to obtain the lower
bound an assumption of uniqueness of equilibrium states is
necessary and this assumption is also used to prove that the
upper bound is strictly negative (see also \cite{Yo90} for
uniformly expanding transformations and for partially
hyperbolic attractors for diffeomorphisms).

Hence the assumption that $\EE$ is formed by a single
measure is natural in this setting. The author feels this
assumption should not be needed to obtain an expression for
the upper bound in terms of entropies, as in \cite{kifer90}.
However the relevant ``volume lemmas'' are presently not
available in the setting of special flows over non-uniformly
expanding base, with singularities or criticalities.
Moreover the uniqueness of equilibrium states with respect
to a large family of potentials (or observables) is still
unknown in general (see
\cite{OliVi2005,AMSV06,arbieto-matheus2006} for recent
progress in this direction). Therefore instead of following
the approach of \cite{kifer90} we have reduced the problem
of estimating the deviations for the suspension flow, with
respect to a continuous observable, to the problem of
estimating deviations for the base transformation, with
respect to an unbounded observable, and then rely on
previous work \cite{araujo-pacifico2006} for non-uniformly
expanding transformations.  To deal with the dynamics near
the singularities we impose conditions of very slow
recurrence to the singular set $\cS$ for the base
transformation $f$ together with a growth condition on the
roof function $r$ near the singularities. In the end to
conclude that the upper bound is strictly negative we use
uniqueness of the relevant equilibrium state. Unfortunately
this argument does not rule out superexponential decay
in~\eqref{eq:flowDeviation}.

Recently Melbourne and Nicol~\cite{melnicol07} obtained
sharp large deviation bounds (i.e. they showed that the
limit \eqref{eq:flowDeviation} exists) for systems modeled
on Markov towers (also known as Young towers) without
requiring uniqueness of equilibrium states. In the same work
upper large deviation bounds are obtained for semiflows over
Markov towers \emph{assuming that the roof function is
  bounded}. However their method presents two disadvantages:
the large deviation estimates in~\cite{melnicol07} are
proved only for H\"older observables, and these estimates are
obtained for the invariant physical measure rather than the
volume or Lebesgue measure, which is more directly
accessible.

Section~\ref{sec:exampl-applic} shows how the conditions of
$f$ and on $r$ are satisfied by many relevant examples. In
particular in Subsection~\ref{sec:singul-hyperb-attrac} it
is explained how to obtain a large deviation bound for
geometric Lorenz flows using the statement of the Main
Theorem applied to suspensions semiflows over piecewise
expanding maps with singularities, which are treated in a
preliminary fashion in
Subsection~\ref{sec:semifl-over-non-1} and at length in
Section~\ref{sec:exponent-slow-approx}. The main result
needed for the proof of the Main Theorem is a large
deviation bound for observables with logarithmic growth near
the singular set for a non-uniformly expanding map, which is
proved in Section~\ref{sec:large-deviat-observ}.
Then the statement of the Main Theorem about large
deviations for a suspension semiflow is reduced to a
statement of large deviations for the dynamics of the base
transformation in Section~\ref{sec:large-deviat-dynamic}
concluding the proof of the Main Theorem.  Note that in
contrast to the results on decay of correlations for Anosov
flows or Axiom A flows, here we do not need any coboundary
conditions on the roof function for the large deviation
bound to hold.  

In Section \ref{sec:isolat-physic-measur} we present a
derivation of the Entropy Formula for non-uniformly
expanding maps with slow recurrence to the singular set,
which is used to establish that some examples presented in
Section~\ref{sec:exampl-applic} do satisfy our assumptions
and which might be interesting in itself.

\subsection*{Acknowledgments}
This paper answers a question posed to me by M. Viana. I am
thankful to S. Senti for having pointed to me how to make
the example in Subsection~\ref{sec:semifl-over-quadrat}
work. The author is indebted to the openness and hospitality
of IMPA, and to its superb library, where most of this work
was written. I wish to thank the referee for his/her
valuable suggestions which helped improve the presentation.


\section{Examples of application}
\label{sec:exampl-applic}

Here we present some concrete examples where our results can
be applied.

\subsection{Suspension semiflows over multidimensional
  volume expanding and quasi-expanding maps}
\label{sec:semifl-over-volume}

Let $f:M\setminus\cS\to M$ be a transitive non-uniformly
expanding map with exponentially slow recurrence to $\cS$
satisfying $J=|\det Df|>1$, $\psi=\log\| (Df)^{-1} \|\le0$
and $\psi=0$ at finitely many points only (a
\emph{quasi-expanding} map). We claim that in this setting
$\EE$ is a singleton.

Indeed $\EE$ is non-empty by Theorem~\ref{thm:abv} since
every absolutely continuous invariant probability measure is
an equilibrium state (see
e.g. Theorem~\ref{thm:reciprocal-EF} in
Section~\ref{sec:isolat-physic-measur}). Since $f\mid
M\setminus \cS$ is a local diffeomorphism and the support of
such absolutely continuous invariant measures contains open
sets, the transitivity together with regularity of the map
ensure that there exists only one absolutely continuous
invariant measure. For otherwise let $\mu_i$ be ergodic
absolutely continuous $f$-invariant probability measures and
let $B_i\subset\supp(\mu_i)$ be open sets in the support
$i=1,2$; by transitivity and continuity there exists a
non-empty open subset $B\subset B_1$ and an iterate such
that $f^n(B)\subset B_2$ and by smoothness $\Leb$-almost
every point in $f^n(B)$ is both a $\mu_1$-generic point and
a $\mu_2$-generic point, thus $\mu_1\equiv\mu_2$. This shows
that there exists a unique absolutely continuous invariant
probability measure for $f$.

Note now that every equilibrium state $\nu\in\EE$ must be
such that $h_\nu(f)=\nu(\log J)>0$ and since $\psi\le0$ and
has at most finitely many zeroes, then either $\nu(\psi)<0$
and by Theorem \ref{thm:reciprocal-EF} the measure $\nu$
must be absolutely continuous, or $\nu(\psi)=0$ and
$\supp\nu\subseteq\psi^{-1}(\{0\})$ is finite thus
$h_\nu(f)=0$, a contradiction.

Therefore by the uniqueness result above $\nu$
must coincide with $\mu$. We have shown that $\EE=\{\mu\}$,
as claimed.

Hence we can apply Theorem~\ref{mthm:LDNUEflow} for
semiflows over non-uniformly expanding maps with
exponentially slow recurrence to the singular set which are
also transitive, volume expanding and expanding except at
finitely many points, and whose roof function grows with the
logarithm of the distance to $\cS$.

For examples of multidimensional local diffeomorphisms in
this setting see \cite{ArTah}. In this case $\cS=\emptyset$
and we can apply Theorem~\ref{mthm:LDNUEflow} for semiflows
with this type of base transformation plus a continuous
(and thus bounded) roof function.

Clearly the same large deviation bound holds for a semiflow
over a local diffeomorphisms which is uniformly expanding
together with any continuous roof function.

\subsection{Suspension semiflows over piecewise
  expanding maps with singularities}
\label{sec:semifl-over-non-1}

Let $M$ be the circle $\sS^1$ or the interval $[0,1]$ with
$\{0,1\}\subset\cS$ and $\cS\subset M$ an at most
denumerable and non-flat singular set of $f$ such that its
closure $\overline{\cS}$ has zero Lebesgue measure:
$\Leb(\overline{\cS})=0$.

If we assume that $-\infty<\psi<-c<0$ on $M\setminus\cS$ for
some $c>0$ (so that in particular there are no critical
points: $\cC=\emptyset$) and that $f$ is transitive with
slow recurrence to $\cS$, then the set $\EE$ of equilibrium
states with respect to $\log |f'|$ is formed by a single
absolutely continuous invariant probability measure, as
shown in Subsection~\ref{sec:semifl-over-volume}, since $f$
is automatically non-uniformly expanding, quasi-expanding
and volume expanding as well.

Observe that for $C^2$ maps in our conditions with finitely
many smoothness domains, or with derivative of bounded
variation, it is well known that there exists a unique
ergodic absolutely continuous invariant probability measure
$\mu$ with bounded density \cite{HK90,Ry83}.  Since the
function $\log\dist(x,\cS)$ is $\Leb$-integrable we also
have that this function is $\mu$-integrable. Thus for all
$\epsilon>0$ there is $\delta>0$ such that
$\int\big|\log\dist_\delta(x,\cS)\big| \, d\mu(x) <
\epsilon.$ By the ergodicity and absolute continuity of
$\mu$ this means that $f$ has slow recurrence to $\cS$ for a
positive Lebesgue measure subset of $M$.
Theorem~\ref{thm:abv} together with \cite{ABV00} ensure that
$f$ is in fact non-uniformly expanding with slow recurrence
to $\cS$. Moreover by \cite{Ke85} the same argument applies
to $C^{1+\alpha}$ piecewise expanding maps with finitely
many smoothness domains, for some $\alpha\in(0,1)$.

To be able to apply the Main Theorem we need exponentially
slow recurrence to $\cS$. We prove this in
Section~\ref{sec:exponent-slow-approx} assuming that $|f'|$
grows as the inverse of some power of the distance to
$\cS'=\cS\cap f(M)$, i.e. besides conditions (S1) through
(S4) we impose
\begin{itemize}
\item[(S5)] $\big|f'(x)\big|\ge
  B^{-1}\dist(x,\cS')^{-\beta}$ for all $x\in
  M\setminus\cS$,
\end{itemize}
where $\cS'$ is the (sub)set of singularities which matters
for the asymptotic dynamics of $f$.


Hence \emph{a semiflow over a piecewise expanding map
  with singularities satisfying some technical conditions,
  and with a roof function having logarithmic growth near
  the singularities} admits a large deviation bound as in
Theorem~\ref{mthm:LDNUEflow}.






\subsection{Suspension semiflows over quadratic maps on
  Benedicks-Carleson parameters}
\label{sec:semifl-over-quadrat}

Set $M=I=[-1,1]$ and suppose the transformation $f$ is given
by $f_a(x)=a-x^2$ for $a\in[a_0,2]$ in the positive Lebesgue
measure subset constructed by Benedicks and Carleson in
\cite{BC85,BC91}, where $a_0\approx2$. The properties of the
family $f_a$ have been thoroughly studied by a considerable
number of people. We just mention that Freitas in
\cite{freitas} showed that for these parameters $f_a$ is not
only a non-uniformly expanding map with $\cS=\cC=\{0\}$ but
also exhibits exponentially slow approximation to the
singular set. Actually in \cite{freitas} only
\emph{subexponentially} slow approximation is stated but the
same arguments yield an exponential bound as well, as
obtained in a much more delicate setting with infinitely
many critical points in \cite{ArPa04}.

Moreover Bruin and Keller \cite{bruin-keller1990} show that
for this class of maps (specifically for
\emph{Collet-Eckman maps}, i.e. such that
$\liminf_{n\to\infty}\big|(f_a^n)'(a)\big|^{1/n}>1$ without
extra conditions of recurrence to the criticality) the
unique absolutely continuous invariant probability measure
is also the unique equilibrium state with respect to $\log
|f_a'|$.

Therefore for any given suspension semiflow over such quadratic
maps $f_a$ with roof function having logarithmic growth near
$0$ we can apply Theorem~\ref{mthm:LDNUEflow}, and obtain a
large deviation bound for these special flows.


\subsection{Lorenz and geometric Lorenz attractors}
\label{sec:singul-hyperb-attrac}

The $C^{1+\alpha}$ map $f:[-1,1]\setminus\{0\}\to[-1,1]$
obtained as the quotient map of the Poincar\'e first return
map $R$ presented in Section~\ref{sec:lorentz-geometr-lore}
through projection along the leaves of the stable foliation
satisfies the following conditions, which define a
\emph{Lorenz-like map}:
\begin{enumerate}
\item there are constants $c>0$ and $\sigma>1$ such that for
  every $n\ge1$ and for all $x\in [-1,1]\setminus\cup_{0\le
    j<n} f^{-n}\{0\}$ we have $\big|(f^n)'(x)\big|\ge
  c\sigma^n$;
\item $f$ has a dense orbit;
\item $f(0^+)=-1$, $f(0^-)=1$, $f(1)\in(0,1)$ and $f(-1)\in(-1,0)$.
\end{enumerate}
Note in particular that there are no critical points and
that for some $k\ge1$ the map $g=f^k$ satisfies the
conditions of Section~\ref{sec:semifl-over-non-1}. (If
$\sigma>\sqrt2$ then $f$ is even locally eventually onto, see
e.g.  \cite{LMP05}, thus transitive.) For exponentially slow
recurrence to the singularities see
Section~\ref{sec:exponent-slow-approx}.  So we can obtain a
large deviation bound for $g$ which easily gives a large
deviation bound for $f$.

Indeed, assume without loss of generality that $\mu(\vfi)=0$
and that for all small $\epsilon>0$ we have
$\Leb\{S_n^g\vfi>n\epsilon\}<C e^{-\zeta n}$ for some
$C(\epsilon),\zeta(\epsilon)>0$ and every $n>0$. It is
enough to argue for a bounded and continuous $\vfi$ as
explained in Section~\ref{sec:large-deviat-observ}.  Then
for $m>0$ we can write $m=nk+p$ with $n>0$ and $0\le p <
k-1$ and also
\begin{align*}
  \frac1m S_m^f\vfi
  &=
  \frac1{nk+p} \big( S^f_p (\vfi\circ f^{nk}) +  S^f_{nk}\vfi
  \big)
  =
  \frac1{nk+p} \big( S^f_p (\vfi\circ f^{nk})
  + \sum_{i=0}^{n-1} S^g_n(\vfi\circ f^i)
  \big)
  \\
  &\le
  \frac{p\sup|\vfi|}{nk+p}
  +
  \frac1{k+p/n} \sum_{i=0}^{n-1} \frac1n  S^g_n(\vfi\circ
  f^i)
  \le
  \frac{p}m \sup|\vfi|
  +
  \frac1k \sum_{i=0}^{n-1} \frac1n  S^g_n(\vfi\circ
  f^i).
\end{align*}
Given $\epsilon>0$ take $m$ so big that
$p\sup|\vfi|/m<\epsilon/2$, note that $\mu(\vfi\circ
f^i)=0$ for all $i\ge0$ and
\begin{align*}
  \big\{ \frac1m S_m^f\vfi >\epsilon \big\}
  \subseteq
  \bigcup_{i=0}^{n-1} \big\{ \frac1n  S^g_n(\vfi\circ
  f^i) > \frac{\epsilon}{2k} \big\}.
\end{align*}
This shows how to reduce the problem of large deviations for
bounded observables to the same problem for a finite power
of the transformation.

To deduce Corollary~\ref{mcor:LDLorenz}, since the reduction
to a large deviation bound for the map $f$ is the content of
Section~\ref{sec:large-deviat-dynamic}, all we need to do
here is to explain how we deduce a large deviation bound for
$R$ from a similar bound for the map $f$. For this we
strongly use the uniform contraction along the leaves of the
stable foliation on the global cross-section $S$ to obtain
the following relation. Denote by $P:S\to[-1,1]$ the
projection $(x,y,1)\mapsto x$.

\begin{lemma}
  \label{le:uniformclose}
  Let $\epsilon>0$ and a bounded continuous function
  $\psi:U\to\RR$ be given in a neighborhood $U$ of the
  geometric Lorenz attractor $\Lambda$. Define
  $\vfi:S\setminus\ell\to\RR$ by $\vfi(x,y,1)=
  \int_0^{\tau(x,y,1)}\psi\big(X^t(x,y,1)\big)\,dt$, where
  $\tau(x,y,1)$ is the first return time to $S$ of the point
  $(x,y,1)\in S$.  Assume without loss of generality that
  $\mu(\vfi)=0$ where $\mu$ is a $R$-invariant probability
  measure such that $\tau$ is $\mu$-integrable.

  Then there exist integers $N,k>1$, a small $\delta>0$, a
  constant $\gamma>0$ dependent on $\psi$ and the flow only,
  and a continuous function
  $l:[-1,1]\setminus\cup_{i=0}^{k-1} f^{-i}\{0\} \to\RR$
  with logarithmic growth near the set
  $\cS_k=\cup_{i=0}^{k-1} f^{-i}\{0\}$ such that for all
  $n>N$
\begin{align}\label{eq:uniformclose}
  \Big\{\big|\frac1n S_n^{R^k}\vfi
  \big|>3\epsilon\Big\}
  \subseteq
    P^{-1}\Big(\Big\{\frac1n
  S_n^{f^k}\Delta_\delta>\frac{\epsilon}{\gamma}
  \Big\}
  \cup
  \Big\{\big|\frac1n S_n^{f^k} l \big|>\epsilon\Big\}\Big).
\end{align}
\end{lemma}

This reduces the problem of estimating the Lebesgue measure
of the left hand side set in~\eqref{eq:uniformclose} to the
estimation of the measure of the right hand side set,
transferring the problem to the dynamics of $g=f^k$, which
is the subject of Section~\ref{sec:semifl-over-non-1} and
Section~\ref{sec:exponent-slow-approx}.

\begin{proof}
  According to the construction of geometric Lorenz flows,
  there are $C>0$ and $0<\lambda<1$ such that given
  $x\in[-1,1]\setminus\{0\}$ and two distinct values
  $y_1,y_2\in[-1,1]$
  \begin{align}\label{eq:contraction}
    \dist\big(R^k(x,y_1,1),R^k(x,y_2,1)\big)\le C \lambda^k
    \text{ for all } 1\le k \le n,
  \end{align}
  where $n\ge1$ is the first time the orbit of the points
  hit the singular line, corresponding to the stable
  foliation of the singularity of the flow. These
  \emph{hitting times} depend only on the orbit of $x$ under
  the map $f$ and correspond to times $n$ for which
  $f^n(x)=0$. But $X_0=\cup_{n\ge0}f^{-n}(\{0\})$ is
  denumerable. Thus the corresponding set of points in $S$,
  given by the lines $\{x\}\times[-1,1]\times\{1\}$ for
  $x\in X_0$, has zero area on $S$. Therefore for a full
  Lebesgue measure subset of $S$ we
  have~\eqref{eq:contraction} for all $k\ge1$.


  Moreover since $(x,y_1,1),(x,y_2,1)$ belong to the same
  stable manifold, then for all times $t>0$ we have
  \begin{align}
    \label{eq:closeWs}
    \dist\big(X^t(x,y_1,1),X^t(x,y_2,1)\big)\le \kappa\cdot
    |y_1-y_2|,
  \end{align}
  for a constant $\kappa>0$ depending only on the angles
  between the surface $S$ and the stable leaves of the flow
  through points of $S$ (which is uniformly bounded by the
  compactness of $S$).  Note that $\vfi$ is continuous on
  $S\setminus\ell$ and
  \begin{align}\label{eq:boundvfi}
    |\vfi(x,y,1)|\le \tau(x,y,1)\cdot\sup|\psi| \le
    -C_0\cdot\log|x|\cdot\sup|\psi|
  \end{align}
  for a constant $C_0>0$, since $\tau$ grows near $\ell$
  like the logarithm of the distance to $\ell$. Then it is
  clear that for $y,w\in[-1,1]$ and $n>1$
  \begin{align*}
    \left|\frac1{n}\sum_{j=0}^{n-1}
      \big(\vfi(R^j(x,y,1))-\vfi(R^j(x,w,1)\big) \right| \le
    \frac1n\sum_{j=0}^{n-1}\overline\vfi_j(x)
  \end{align*}
  where $(x_j,y_j,1)=R^j(x,y,1)$ for $j\ge0$, $(x,y,1)\in
  S$, and 
  \begin{align*}
    \overline\vfi_j(x)&=\sup_{y,w\in[-1,1]}
    \big|\vfi\big(R^j(x,y,1)\big)-\vfi\big(R^j(x,w,1)\big)\big|.
  \end{align*}
  
  Let $\epsilon>0$ be given. Choose a small $\delta>0$ and
  $\eta>0$ such that $-C_0\kappa\eta\log\delta<\epsilon/3$
  and $\kappa\eta\le\sup|\psi|$. Let $\xi>0$ satisfy
  \begin{align}\label{eq:unifcont}
    \dist\big( (x,y,z) , (x',y',z') \big) < \xi \implies |
    \psi(x,y,z) - \psi(x',y',z') | <\eta.
  \end{align}
  Then we may find by \eqref{eq:contraction} a
  $j_0=j_0(\eta)\ge1$ such that $|y_j-w_j|\le\xi/\kappa$ for
  $j>j_0$ and any pair $y,w$ in the same vertical line.
  Thus we also get
  after~\eqref{eq:closeWs},~\eqref{eq:unifcont} and the
  choices of $\epsilon,\delta$ and $\eta$
  \begin{align}\label{eq:overline}
    \overline\vfi_j(x) 
    &\le
    -C_0\log|x_j|\cdot\sup_{0<t<-C_0\log|x_j|}
    \big|\psi\big(X^t(x_j,y_j,1)\big)-\psi\big(X^t(x_j,w_j,1)\big)\big|
    \nonumber
    \\
    &\le -C_0\log|x_j|\cdot \kappa\eta 
    \le
    C_0 \kappa\eta\cdot \Delta_\delta(x_j) +\epsilon/2.
  \end{align}
  Take a continuous $l:[-1,1]\setminus\{0\}\to\RR$ such that
  for some $0<a<\epsilon/3$
  \begin{enumerate}
  \item $\min_{y\in[-1,1]}\vfi\big(R^{j_0}(x,y,1)\big)-a \le
    l(x)\le
    a+\max_{y\in[-1,1]}\vfi\big(R^{j_0}(x,y,1)\big)$; and
  \item $\mu(l\circ P)=\mu(\vfi)$.
  \end{enumerate}
  Note that $\vfi$ is $\mu$-integrable: this follows
  from the boundedness assumption on $\psi$ and by the
  $\mu$-integrability of $\tau$ after~\eqref{eq:boundvfi}.
  Observe that $l$ as above has logarithmic growth near
  $\cS_k$ by definition.

  To obtain such function $l$ disintegrate $\mu$ along the
  measurable partition of $S$ given by the vertical lines
  $\{x\}\times[-1,1]\times\{1\}$ and define $l_0(x)=\int
  \vfi\,d\mu_x$. Then approximate $l_0$ by a continuous
  function $l_1$ such that $\int |l_0-l_1|\circ P
  \,d\mu<\epsilon/3$ (through e.g. a convolution).  Now for
  some $-\epsilon/3<a<\epsilon/3$ the function $l=l_1+a$
  satisfies conditions 1-2 above.

  Now for $n>0$ using~\eqref{eq:overline}, $f\circ P=P\circ
  R$ and summing over orbits of $R^k$ and $f^k$
  \begin{align}
    |S_n(l\circ P)-S_n\vfi|(x,y,1)
    &\le 
    | l\circ P - \vfi|(x,y,1) + |S_{n-1} (l\circ P -
    \vfi)|(x,y,1)\nonumber
    \\
    &\le 
    2\sup|\psi| C_0\log|x| +a+
    \sum_{j=1}^{n-1}
    \big(C_0\kappa\eta\Delta_\delta(f^{jk}(x)) + \frac\epsilon3
    +a\big)
    \nonumber
    \\
    &\le 
    2\sup|\psi| C_0 (\log \delta^{-1} +
    \Delta_\delta(x))+ \frac{2n\epsilon}3 +
    C_0\kappa\eta \cdot S_{n-1}\Delta_\delta(f^k (x)) \nonumber
    \\
    &\le
    2\sup|\psi| C_0 \log \delta^{-1} + \frac{2n\epsilon}3 +
    C_0(\kappa\eta +2\sup|\psi|)\cdot S_{n}\Delta_\delta(x). \label{eq:zeta1}
  \end{align}
  Observe  that
  \begin{align}\label{eq:zeta2}
    \big\{ \big|\frac1n S_n^{R^k}\vfi \big|>3\epsilon\big\}
    \subseteq
    \big\{ \big| \frac1n\big( S_n^{R^k}(l\circ P)-S_n^{R^k}\vfi\big)
    \big| > 2\epsilon \}
    \cup
    \big\{ \frac1n\big| S_n^{R^k}(l\circ P)\big|>\epsilon\big\}.
  \end{align}
  From~\eqref{eq:zeta1}, setting $\gamma_1=2\sup|\psi| C_0
  \log\delta^{-1}$ and $\gamma_2=C_0(\kappa\eta
  +2\sup|\psi|)$ we obtain for $n$ big enough
  \begin{align*}
    \frac1n\big( S_n^{R^k}(l\circ P)-S_n^{R^k}\vfi\big)
    \le
    \frac{2\epsilon}3 + \frac{\gamma_1}n +
    \frac{\gamma_2}n\cdot S_{n}^{f^k}\Delta_\delta\circ P
    \le
    \epsilon + \frac{\gamma_2}n\cdot
    S_{n}^{f^k}\Delta_\delta \circ P
  \end{align*}
  where $\gamma_2\le 3C_0\sup|\psi|$ by the choice of
  $\eta$. Hence
  \begin{align*}
    \big\{ \big| \frac1n\big( S_n^{R^k}(l\circ P)-S_n^{R^k}\vfi\big)
    \big| > 2\epsilon \}
    \subseteq P^{-1}
    \big\{ \frac1n S_{n}^{f^k}\Delta_\delta > 
    \frac{\epsilon}{3C_0\sup|\psi|} \big\}
  \end{align*}
  and this together with~\eqref{eq:zeta2} completes the
  proof of the lemma.
\end{proof}



\section{Large deviations for observables with logarithmic
  growth near singularities}
\label{sec:large-deviat-observ}

The main bound on large deviations for suspension semiflows
over a non-uniformly expanding base will be obtained from
the following large deviation statement for non-uniformly
expanding transformations. 

\begin{theorem}
  \label{thm:LDNUElog}
  Let $f:M\to M$ be a regular $C^{1+\alpha}$ local
  diffeomorphism on $M\setminus \cS$ where $\cS$ is a
  non-flat critical set and $\alpha\in(0,1)$.
  Assume that $f$ is a non-uniformly expanding map with
  exponentially slow recurrence to the singular set $\cS$
  and let $\vfi:M\setminus\cS\to\RR$ be a continuous map
  which has logarithmic growth near $\cS$. Moreover assume
  that there exists a unique equilibrium state $\mu$ with
  respect to $\log J$ which is absolutely continuous. Then
  for any given $\omega>0$
  \begin{align*}
    \limsup_{n\to+\infty}\frac1n\log\Leb\Big\{ x\in M:
    \left|\frac1n S_n\vfi(x)-\mu(\vfi)\right|\ge\omega
    \Big\}<0.
\end{align*}
\end{theorem}

\begin{proof}
Define
\begin{align*}
 \vfi_k=\xi_k\circ\vfi
\mbox{  where  }
\xi_k(x)=\left\{
  \begin{array}[l]{ll}
k & \mbox{if  } x\ge k
\\
x & \mbox{if  } |x|<k
\\
-k & \mbox{if  } x\le -k
  \end{array}
\right.,\,\, k\ge1.
\end{align*}
Then $\vfi_k:M\to\RR$ is continuous for all $k\ge1$,
$\vfi_k(x)\xrightarrow[k\to\infty]{}\vfi(x)$ for all $x\in
M\setminus\cS$ and
$|\vfi-\vfi_k|\le\big|\vfi\big|\chi_{\{|\vfi|>k\}}$. Moreover
we clearly have for all $n,k\ge1$
\begin{align}
  \label{eq:approxlog}
S_n\vfi_k- S_n\big|\vfi-\vfi_k\big| \le
S_n\vfi=S_n\vfi_k+S_n(\vfi-\vfi_k) \le S_n\vfi_k +
S_n\big|\vfi-\vfi_k\big|.
\end{align}
Observe that, since $\vfi$ has logarithmic growth near $\cS$
(see \eqref{eq:log-growth}), for any given $c,\epsilon_0>0$
we may choose $\epsilon_1,\delta_1>0$ such that the
exponential slow recurrence condition
\eqref{eq:expslowrecurrence} is true and
$K\cdot\epsilon_1\le\epsilon_0$. Then choose $k\ge1$ very
big so that $\{|\vfi|>k\}\subseteq B(\cS,\delta_1)$.  From
\eqref{eq:approxlog} we obtain the following inclusions
\begin{align}
  \left\{\frac1n S_n\vfi > c\right\}
  &\subseteq
  \left\{\frac1n S_n\vfi_k + \frac1n S_n\big|\vfi-\vfi_k\big| >
    c\right\}
  \subseteq
  \left\{\frac1n S_n\vfi_k > c- K\epsilon_1\right\}
  \cup
  \left\{\frac1n S_n\big|\vfi-\vfi_k\big| > K\epsilon_1\right\}\nonumber
  \\
  &\subseteq
  \left\{\frac1n S_n\vfi_k > c- \epsilon_0\right\}
  \cup
  \left\{\frac1n S_n\Delta_{\delta_1}\ge\epsilon_1\right\},
  \label{eq:approxlog1}
\end{align}
where in \eqref{eq:approxlog1} we use the assumption that
$\vfi$ is of logarithmic growth near $\cS$ and the choices
of $\epsilon_1,\delta_1$. Analogously we get with opposite
inequalities
\begin{align}
  \left\{\frac1n S_n\vfi < c \right\}
  &\subseteq
  \left\{\frac1n S_n\vfi_k - \frac1n S_n\big|\vfi-\vfi_k\big| <
    c\right\}
  \subseteq
  \left\{\frac1n S_n\vfi_k < c + K\epsilon_1\right\}
  \cup
  \left\{\frac1n S_n\big|\vfi-\vfi_k\big|> K\epsilon_1\right\}\nonumber
  \\
  &\subseteq
  \left\{\frac1n S_n\vfi_k < c + \epsilon_0\right\}
  \cup
  \left\{\frac1n S_n\Delta_{\delta_1}\ge\epsilon_1\right\}.
  \label{eq:approxlog2}
\end{align}
From \eqref{eq:approxlog1} and \eqref{eq:approxlog2} we see
that \emph{to get the bound for large deviations in the statement
of Theorem~\ref{thm:LDNUElog} it suffices to obtain a large
deviation bound for the continuous function $\vfi_k$ with
respect to the same transformation $f$} and \emph{to have
exponentially slow recurrence to the singular set $\cS$}.

To obtain this large deviation bound, we use the following
result already obtained for continuous observables over
non-uniformly expanding transformations in our setting, see
\cite{araujo-pacifico2006}.
\begin{theorem}
  \label{thm:LDNUE}
  Let $f:M\to M$ be a local diffeomorphism outside a
  non-flat singular set $\cS$ which is non-uniformly
  expanding and has slow recurrence to $\cS$.  For
  $\omega_0>0$ and a
  continuous function $\vfi_0:M\to\RR$ 
  there exists $\epsilon,\delta>0$ arbitrarily close to $0$
  such that, writing
\begin{align*}
A_n=\{x\in M: \frac1nS_n\Delta_\delta(x)\le\epsilon\}
\quad\text{and}\quad
B_n=\left\{
x\in M : 
\inf\big\{\big|
\frac1n S_n\vfi_0(x) - \eta(\vfi_0)
\big| : \eta\in\EE \big\}
> \omega_0
\right\}
\end{align*}
we get $
\limsup_{n\to+\infty}\frac1n
\log \Leb\big(A_n\cap B_n\big) <0.$
\end{theorem}
Recall that $\EE$
is the set of all equilibrium states of $f$ with respect to
the potential $\log J$.

Note that exponentially slow recurrence implies
$\limsup_{n\to+\infty}\frac1n\Leb(M\setminus A_n) <0$. Under
this assumption Theorem~\ref{thm:LDNUE} ensures that for
$(\epsilon,\delta)$ close enough to $(0,0)$ we get
$\limsup_{n\to+\infty}\frac1n \log \Leb( B_n) <0.$ To use
this we also need that \emph{$\EE$ consists only of the
  unique absolutely continuous invariant probability measure
  $\mu$}.
Under this uniqueness assumption we have $\EE=\{\mu\}$ in
Theorem~\ref{thm:LDNUE} and take $\omega,\epsilon_0>0$
small, choose $k\ge1$ as before, set $\vfi_0=\vfi_k$ and
$\omega_0=\omega+\epsilon_0$.  In \eqref{eq:approxlog1} set
$c=\mu(\vfi_0)-\omega$ and in \eqref{eq:approxlog2} set
$c=\mu(\vfi_0)+\omega$. Then we have the inclusion
\begin{align}\label{eq:omegazero}
  \left\{\Big|\frac1n
    S_n\vfi-\mu(\vfi)\Big|>\omega\right\}
  \subseteq
  \left\{\Big|\frac1n
    S_n\vfi_0-\mu(\vfi_0)\Big|>\omega_0\right\}
   \cup
  \left\{\frac1n S_n\Delta_{\delta_1}\ge\epsilon_1\right\}.
\end{align}
By Theorem~\ref{thm:LDNUE} we may find $\epsilon,\delta>0$
small enough so that the exponentially slow recurrence holds
also for the pair $(\epsilon,\delta)$ and hence
\begin{align}
  \label{eq:omegazero1}
\limsup_{n\to+\infty}\frac1n\log\Leb
\left\{\Big|\frac1n
S_n\vfi_0-\mu(\vfi_0)\Big|>\omega_0\right\}<0.
\end{align}
Finally the choice of $\epsilon_1,\delta_1$ according to the
condition on exponential slow recurrence to $\cS$ ensures
that the Lebesgue measure of the right hand subset in
\eqref{eq:omegazero} is also exponentially small when
$n\to\infty$. This together with \eqref{eq:omegazero1}
concludes the proof of Theorem~\ref{thm:LDNUElog}.
\end{proof}



\section{Large deviations and the dynamics on the base}
\label{sec:large-deviat-dynamic}

Here we show how the large deviation bound for a semiflow
over a non-uniformly expanding base can be deduced from a
large deviation bound for the base dynamics, under a
logarithmic growth condition on the roof function.

\subsection{Reduction to the base dynamics}
\label{sec:reduct-argument-base}

Let $\psi:M_r\to\RR$ be continuous and bounded.  For $T>0$
and $z=(x,s)$ with $x\in M$ and $0\le s < r(x) <\infty$ we
can write
\begin{align*}
\int_0^T
\hspace{-0.2cm}\psi\big(X^t(z)\big)\,dt
=
\int_s^{r(x)}
\hspace{-0.6cm}\psi\big(X^t(x,0)\big)\,dt
+
\sum_{j=1}^{n-1} \int_0^{r(f^j(x))} 
\hspace{-0.6cm}\psi\big(X^t(f^j(x),0)\big)\,dt
+\int_0^{T+s-S_{n}r(x)}
\hspace{-0.6cm}\psi\big(X^t(f^n(x),0)\big)\,dt,
\end{align*}
where $n=n(x,s,T)\in\NN$ is the ``lap number''  such that $S_{n}r(x)\le s+T <
S_{n+1} r(x)$.

Setting $\vfi(x)=\int_0^{r(x)}\psi(x,0)\,dt$ we obtain
\begin{align*}
  \frac1T\int_0^T \hspace{-0.2cm}\psi\big(X^t(z)\big)\,dt
  =\frac1T S_n\vfi(x)
  -\frac1T\int_0^s\psi\big(X^t(x,0)\big)\,dt
  +\frac1T\int_0^{T+s-S_n r(x)}
  \hspace{-0.6cm}\psi\big(X^t(f^n(x),0)\big)\,dt.
\end{align*}
Clearly we can bound the sum $I=I(x,s,T)$ of the two
integral terms on the right hand side above by
\begin{align}\label{eq:resto}
  I= I(x,s,T)\le \left(2\frac{s}T + \frac{S_{n+1}r(x)-S_n
      r(x)}T\right)\cdot\|\psi\|,
\end{align}
where $\|\psi\|=\sup|\psi|$.  Observe that for a given
$\omega>0$ and for $0\le s<r(x)$ and $n=n(x,s,T)$
\begin{align}\label{eq:LD0}
  \left\{ (x,s)\in M_r : \left| \frac1T S_n\vfi(x) + I(x,s,T)
    -\frac{\mu(\vfi)}{\mu(r)}\right| > \omega\right\}
\end{align}
is contained in 
\begin{align}\label{eq:LD1}
  \left\{ (x,s)\in M_r : \left| \frac1T S_n\vfi(x)
      -\frac{\mu(\vfi)}{\mu(r)}\right| > \frac\omega2 \right\}
\cup
\left\{ (x,s)\in M_r: I(x,s,T) >\frac\omega2\right\}.
\end{align}
Note that if $\psi\equiv0$ then we need only consider the
left hand subset of \eqref{eq:LD1} in what follows.  Now we
bound the $\lambda$-measure of each subset above assuming
that $\psi$ is not identically zero.

\subsection{Using the roof function as an observable over the base dynamics}
\label{sec:using-roof-functi}

We start with the right hand subset in \eqref{eq:LD1}. Take
$N\ge1$ big enough so that $N\|\psi\|>2$ and note that for
any given $T,\omega>0$ using \eqref{eq:resto} and
$n=n(x,s,T)$
\begin{align}
  &\lambda\Big\{I>\frac\omega2\Big\}
  =
  \int d\Leb(x)\int_0^{r(x)} \!\!\! ds\,
  \big(\chi_{(\omega/2,+\infty)}\circ I\big)(x,s,T)\nonumber
  \\
  &\le
  \Leb\left\{ r>\frac{\omega T}{2N\|\psi\|} \right\}
  +
  \frac{\omega T}{2N\|\psi\|}\sum_{i=0}^{[T/r_0]+1}
  \Leb\left\{
    \frac{|S_{i+1}r-S_i r|}T>
    \frac{N\|\psi\|-2}{2N\|\psi\|}\omega
    \right\},\label{eq:LD2}
\end{align}
where in the right hand summand we restrict to points $x\in
M$ such that $2N\|\psi\| r(x)\le\omega T$ and $S_ir(x)\le T
< S_{i+1}r(x)$ for each possible lap number $i\in\NN$.
Note that since $r$ is bounded from below $r\ge r_0 >0$ we
have $T\ge r_0 n$ which gives an upper bound $[T/r_0]+1$ for
the possible lap numbers appearing in the summation above,
where $[t]$ denotes $\max\{k\in\ZZ:k\le t\}$, the integer
part of $t\ge0$.  In \eqref{eq:LD2} we have also used the
relations
\begin{align*}
  \frac{2s}T<\frac{2r}T\le\frac{\omega}{N\|\psi\|}
  \quad\text{and}\quad
  \frac{\omega}2-\frac{\omega}{N\|\psi\|}=
  \frac{N\|\psi\|-2}{2N\|\psi\|}\cdot\omega.
\end{align*}
On the one hand, since $r$ grows as the logarithm of the
distance to $\cS$, we have that the left hand summand in
\eqref{eq:LD2} is bounded by
\begin{align}\label{eq:volumebound}
\Leb\left\{x\in M:
\dist\big(x,\cS\big)\le \exp\big(-C\cdot\frac{\omega
T}{2N\|\psi\|}\big)\right\}\le
e^{-C\cdot\kappa\cdot\omega T\|/(2N\|\psi\|)},
\end{align}
where $C>0$ is a constant depending on $r$ only, and we use
condition (S4) on the geometry of $\cS$.  On the other hand,
from $T\ge S_i r(x)\ge r_0 i$ we get the following upper
bound for the summands in the right hand side of
\eqref{eq:LD2} for each $i=0,\dots,[T/r_0]+1$
\begin{align}
  &\Leb\left\{
    \frac{|S_{i+1}r-S_i
      r|}i>\left(\frac{N\|\psi\|-2}{2N\|\psi\|}r_0\right)\cdot\omega 
    \right\} \quad \Big(\text{let }
    r_0'=\frac{N\|\psi\|-2}{2N\|\psi\|}r_0 \Big)\nonumber
    \\
    &\le
    \Leb\left\{\Big|\frac1i S_i r-\mu(r)\Big|>
      \frac{\omega r_0'}2\right\}
    +
    \Leb\left\{\Big|\frac1i S_{i+1}r-\mu(r)\Big|>
      \frac{\omega r_0'}2\right\}
    \le
    2 C_0 e^{-\beta i}\label{eq:LDS11}
\end{align}
for some constants $C_0,\beta>0$, since we have a large
deviation bound for the observable $r$ with respect to the
unique physical measure $\mu$ for $f$.  Recall (see
Section~\ref{sec:large-deviat-observ}) that we took $r$ to
be $\mu$-integrable, continuous on $M\setminus\cS$ and with
logarithmic growth near $\cS$, and $f$ is a non-uniformly
expanding map with exponentially slow recurrence
to 
$\cS$. Consequently we can bound the summation
in~\eqref{eq:LD2} as
\begin{align}
  \frac{\omega T}{2N\|\psi\|}\cdot 2C_0 \sum_{i=0}^{[T/r_0]+1}
  e^{-\beta i}
  \le
  \frac{C \omega T}{2N\|\psi\|} \cdot e^{-\beta T/r_0}\label{eq:LDS12}
\end{align}
for a constant $C>0$ depending on $f,r,\omega$ and $\psi$.
Altogether we see that $\lambda\{I>\omega/2\}$ is bounded by
twice the maximum of the summands in \eqref{eq:LD2}.

From this we obtain
\begin{align}\label{eq:limsup1}
  \limsup_{T\to\infty}
  \frac1T\log\lambda\Big\{I>\frac\omega2\Big\}<0,
\end{align}
as long as we take $\omega>0$ small enough. 

\subsection{Using $\vfi$ as an observable over the base dynamics}
\label{sec:using-vfi-as}

Now for the left hand subset in \eqref{eq:LD1}, note first
that for $\mu$- and $\Leb$-almost every $x\in M$ and every
$0\le s<r(x)$
\begin{align}\label{eq:LD3}
 \frac{S_nr(x)}n\le
 \frac{T+s}{n}\le\frac{S_{n+1}r(x)}n
 \quad\text{so}\quad
 \frac{n(x,s,T)}T\xrightarrow[T\to\infty]{}\frac1{\mu(r)}.
\end{align}
We also have (recall that $n=n(x,s,T)$)
\begin{align*}
  \left| \frac1T S_n\vfi-\frac{\mu(\vfi)}{\mu(r)} \right|
  \le
  \left|\frac{n}T\cdot\frac{S_n\vfi}n-\frac{n}T
    \mu(\vfi)\right|
  +
  \left|\frac{n}T \mu(\vfi) - \frac{\mu(\vfi)}{\mu(r)}
  \right|
  \le
  \frac{n}T\left| \frac{S_n\vfi}n -\mu(\vfi)\right|
  +
  |\mu(\vfi)|\left| \frac{n}T -\frac1{\mu(r)}  \right|.
\end{align*}
Hence the left hand subset in \eqref{eq:LD1} is contained in
\begin{align}
  \label{eq:LD4}
  \left\{(x,s)\in M_r: \frac{n}T\left| \frac{S_n\vfi}n
      -\mu(\vfi)\right|>\frac{\omega}4\right\}
  \cup
  \left\{(x,s)\in M_r: \left| \frac{n}T -\frac1{\mu(r)}
    \right| > \frac{\omega}{4|\mu(\vfi)|}\right\}.
\end{align}
Notice that the $\lambda$-measure of the right hand subset of
\eqref{eq:LD4} is bounded from above by
\begin{align}
  & \lambda\left\{ \left| \frac{n}T -\frac1{\mu(r)} \right|
    > \frac{\omega}{4|\mu(\vfi)|}\,\&\, r\le T\right\} +
  \lambda\{ r>  T \}\nonumber
  \\
  &\le T \sum_{i=0}^{[T/r_0]+1}\sum_{j=0,1} \Leb\left\{x\in M: \left|
      \frac{i}{S_{i+j}r} -\frac1{\mu(r)} \right| >
    \frac{\omega}{|\mu(\vfi)|}\right\} + \int_{\{r> T\}} 
  \hspace{-0.5cm} r  \,\, d\Leb\label{eq:LD41}
\end{align}
where we have used the relation \eqref{eq:LD3}, for small
enough $\omega>0$ and big enough $T$ and $n$. The first
summand in \eqref{eq:LD41} can be bounded using the large
deviation bound for the observable $r$ as before: there are
constants $C_0,\beta>0$ such that
\begin{align*}
   \Leb\left\{x\in M: \left|
      \frac{i}{S_{i+j}r} -\frac1{\mu(r)} \right| >
    \frac{\omega}{|\mu(\vfi)|}\right\} \le C_0e^{-\beta i}
  \quad\text{for}\quad j=0,1,
\end{align*}
and so for some constant $C_1>0$ depending only on $f,r$ and
$\omega$ we get
\begin{align*}
  T \sum_{i=0}^{[T/r_0]+1}\sum_{j=0,1} \Leb\left\{x\in M: \left|
      \frac{i}{S_{i+j}r} -\frac1{\mu(r)} \right| >
    \frac{\omega}{|\mu(\vfi)|}\right\}
  \le C_1 T  e^{-\beta T / r_0}.
\end{align*}
The second summand in \eqref{eq:LD41} is easily bounded
using condition (S4) as follows: for big enough $T>0$ such
that $i>[T]$ implies $(i+1) e^{c_0 i}<1$, where
$c_0=-\kappa\log\rho/(2K)>0$, we have as in
\eqref{eq:volumebound}
\begin{align}
  \int_{\{r> T\}} \hspace{-0.5cm} r  \,\, d\Leb
  &\le
  \sum_{i\ge[T]}\int_i^{i+1} r \,\,d\Leb
  \le
  \sum_{i\ge[T]}(i+1)\Leb\{r>i\}\nonumber
  \\
  &\le
  C_\kappa \sum_{i\ge[T]} (i+1)e^{-2c_0 i}
  \le
  C_\kappa\sum_{i\ge[T]} e^{-c_0 i}
  \le
  C_2\cdot e^{-c_0 T}\label{eq:LD42}
\end{align}
for a positive constant $C_2>0$ depending only on $f$.

Finally the left hand subset of \eqref{eq:LD4} is contained
in the following union
\begin{align}\label{eq:LD5}
  \left\{ (x,s)\in M_r:\left| \frac{T}n
     -\mu(r)\right|>\frac{\mu(r)}2\cdot \omega \right\}
\cup
 \left\{ (x,s)\in M_r: \left|\frac{S_n\vfi}n
      -\mu(\vfi)\right|>\frac{\mu(r)}2\cdot \frac{\omega}4\right\}.
\end{align}
Again for small $\omega>0$ the $\lambda$-measure of the left
hand subset in~\eqref{eq:LD5} is exponentially small with
$T$, using similar arguments to~\eqref{eq:LD41}
and~\eqref{eq:LD42}.  For the right hand subset in
\eqref{eq:LD5} we use the large deviation bound for the
observable $\vfi$ with respect to $f$, since $\vfi$ has also
logarithmic growth near $\cS$. In fact
$\big|\vfi(x)\big|\le\int_0^{r(x)}|\psi(x,s)|\,dt \le
r(x)\cdot\|\psi\| $ for $x\in M\setminus\cS$ because
$\psi:M_r\to\RR$ is bounded.  We can estimate the
$\lambda$-measure of the right hand subset in \eqref{eq:LD5}
as in~\eqref{eq:LD41} through \eqref{eq:LD42} (or as
in~\eqref{eq:LDS11} and \eqref{eq:LDS12}), obtaining
constants $C_3,\gamma>0$ depending on $f,r$ and $\omega$ such that
\begin{align*}
  \lambda\left\{ \left|\frac{S_n\vfi}n
      -\mu(\vfi)\right|>\frac{\mu(r)}2\cdot
    \frac{\omega}4\right\}
  \le C_3 T e^{-\gamma n}.
\end{align*}
From this we conclude
\begin{align}\label{eq:limsup2}
  \limsup_{T\to\infty}
  \frac1T\log\lambda\left\{ \left| \frac1T S_n\vfi
      -\frac{\mu(\vfi)}{\mu(r)}\right| > \frac\omega2
  \right\} < 0.
\end{align}
Putting \eqref{eq:limsup1} and \eqref{eq:limsup2} together,
\emph{as long as we have a result on large deviations for
  continuous observables in $M\setminus\cS$ with logarithmic
  growth near $\cS$, with respect to the dynamics of $f$ and
  the Lebesgue measure}, and \emph{the volume of
  neighborhoods of $\cS$ is comparable to a power of the
  radius}, we are able to prove the Main Theorem for the
suspension flow over $f$.

We have obtained the large deviation bounds needed for the
base dynamics in Section~\ref{sec:large-deviat-observ}, so
the proof of Theorem~\ref{mthm:LDNUEflow} is complete.



\section{The Entropy Formula for non-uniformly expanding maps}
\label{sec:isolat-physic-measur}

Here we obtain the Entropy Formula when $f$ is a
non-uniformly expanding map with slow recurrence to the
singular set. The singular set $\cS$ is formed by critical
points of $f$ and points where $f$ is either not defined, is
not continuous or is not differentiable. Recall from the
Introduction that $\psi=\log\| (Df)^{-1} \|$ and that
$J=|\det Df|$.

\begin{theorem}
  \label{thm:reciprocal-EF}
  Let $f:M\to M$ be a non-uniformly expanding map with slow
  recurrence to the non-flat singular set $\cS$. Let
  $\mu\in\M_f$ be such that $\mu$ is $f$-ergodic,
  $h_\mu(f)=\mu(\log J)$, $-\infty<\mu(\psi)<0$ and for
  every given $\epsilon>0$ there exists $\delta>0$ so that
  $\mu(\Delta_\delta)<\epsilon$. Then $\mu\ll\Leb$ and
  consequently $\mu\in\overline\co(\FF)$.

  Reciprocally, let $\mu\in\M_f$ be such that $\mu$ is
  absolutely continuous with respect to $\Leb$ and assume
  that $\Delta_\delta$ is $\mu$-integrable. Then
  $h_\mu(f)=\mu(\log J)$.
\end{theorem}

Here $\overline{\co}(\FF)$ is the weak$^*$ closure of the
convex hull of the finite set $\FF$ of ergodic physical
probability measures for $f$.  Clearly this is a particular
case of the more general Entropy Formula obtained by
Ledrappier and Young \cite{LY85,LY85a} applied to maps with
singularities and/or criticalities. For $C^2$ endomorphisms
(i.e. smooth maps with criticalities but no singularities)
see Bahnm\"uller and Liu \cite{Li98,BhLi98} for a general
statement. A similar result for piecewise smooth
one-dimensional maps with finitely many branches was
obtained by Ledrappier \cite{Le81}.

As an easy corollary we deduce that $\overline{\co}(\FF)$ is
isolated among the set $\EE$ of all equilibrium states of
$f$ with respect to $J=\log\big|\det Df\big|$, which might
be of independent interest for the ergodic theory of
non-uniformly expanding transformations.

\begin{corollary}
  \label{cor:isolation}
  Let $f:M\to M$ be a non-uniformly expanding map with slow
  recurrence to the non-flat singular set $\cS$. Then there
  exists a weak$^*$ neighborhood $\U$ of $\overline\co(\FF)$ in
  $\M_f$ such that $\U\cap\EE=\overline\co(\FF)$.
\end{corollary}

\begin{proof}
  Take any weak$^*$ neighborhood $\U$ of $\overline\co(\FF)$
  such that every $\mu\in\U$ satisfies $\mu(\psi)<0$. Hence
  every $\mu\in\U\cap\EE$ satisfies the conditions of
  Theorem~\ref{thm:reciprocal-EF}, thus $\mu\in\overline\co(\FF)$.
\end{proof}

Note that whenever the Entropy Formula and its reciprocal
hold for measures close to $\FF$ then the argument proving
Corollary~\ref{cor:isolation} is applicable and we deduce
that $\FF$ is isolated in $\EE$.
The proof of Theorem~\ref{thm:reciprocal-EF} is longer and
occupies the rest of this section. 

\subsection{Hyperbolic times}
\label{sec:hyperb-times}

Here we present some technical results for the study of
non-uniformly expanding maps whose proof can be found in
\cite{Pl72,ABV00,Ze03}.

We say that $n$ is a \emph{$(\sigma,\delta,b)$-hyperbolic
  time of $f$} for a point $x$ if there are $0<\sigma<1$ and
$b,\delta>0$ such that $\prod_{j=n-k}^{n-1}\big\|
Df\big(f^j(x)\big)^{-1}\big\| \le \sigma^k$ and $d_\delta
\big( f^k(x),\cS \big) \ge e^{-bk}$ hold for all
$k=0,\dots,n-1$.

We now outline the properties of these special times. For detailed
proofs see~\cite[Proposition 2.8]{ABV00} and~\cite[Proposition
2.6, Corollary 2.7, Proposition 5.2]{AA03}.

\begin{proposition}
  \label{pr:prophyptimes}
  There are constants $C_1,\delta_1>0$ depending on
  $(\sigma,\de,b)$ and $f$ only such that, if $n$ is
  $(\sigma,\de,b)$-hyperbolic time of $f$ for $x$, then
  there are \emph{hyperbolic pre-balls} $V_k(x)$ which are
  neighborhoods of $f^{n-k}(x)$, $k=1,\dots, n$, satisfying
\begin{enumerate} 
\item $f^k\mid V_k(x)$ maps $V_k(x)$ diffeomorphically to
  the ball of radius $\delta_1$ around $f^n(x)$;
\item 
$d\big(f^{n-k}(y),f^{n-k}(z)\big)\le
  \sigma^{k/2}\cdot d\big(f^{n}(y),f^{n}(z)\big)$
for every $1\leq k\leq n$ and $y,z\in V_k(x)$;
\item 
$C_1^{-1}\le
\big|\det Df^{n-k}(y)\big|/\big|\det Df^{n-k}(z)\big| \le
C_1$ for $y,z\in V_k(x)$.
\end{enumerate}
\end{proposition}

The following ensures existence of infinitely many
hyperbolic times for $\mu$-almost every point for
non-uniformly expanding maps with respect to an ergodic
invariant probability measure $\mu$. A complete proof can be
found in~\cite[Section 5]{ABV00}.

\begin{theorem}
\label{thm:tempos-hip-existem}
Let $f:M\to M$ be a $C^{1+\alpha}$ local diffeomorphism away
from a non-flat singular set $\cS$, for some
$\alpha\in(0,1)$, non-uniformly expanding and with slow
recurrence to $\cS$, with respect to an ergodic invariant
probability measure $\mu$. That is there exists $c>0$ such
that
\begin{align*}
      \limsup_{n\to+\infty}\frac1n S_n\psi \le -c
      \quad\mu-\text{almost everywhere}
\end{align*}
and for every $\epsilon>0$ there exists $\delta>0$ such that
\begin{align*}
    \limsup_{n\to\infty} \frac1n S_n\Delta_\delta(x)\le
    \epsilon
          \quad\mu-\text{almost everywhere}.
\end{align*}
Then there are $\sigma\in(0,1)$, $\delta,b>0$ and there
exists $\theta=\theta(\sigma,\delta,b)>0$ such that
$\mu$-a.e.  $x\in M$ has infinitely many
$(\sigma,\de,b)$-hyperbolic times.  Moreover if we write
$0<n_1<n_2<n_2<\dots$ for the hyperbolic times of $x$ then
their asymptotic frequency satisfies
$\liminf_{N\to\infty}\frac{\#\{ k\ge1 : n_k\le
  N\}}{N}\ge\theta$ for $\Leb\mbox{-a.e.  } x\in M$.
\end{theorem}


\subsection{Existence of generating partition}
\label{sec:existence-generat-pa}

Let $\mu$ be an $f$-invariant ergodic probability measure in
the conditions of the first part of the statement of
Theorem~\ref{thm:reciprocal-EF}.  

Observe first that since $\mu(\psi)<0$ and $\mu$ is ergodic,
then $f$ is non-uniformly expanding. Moreover by the
assumptions on $\mu(\Delta_\delta)$ we see that $f$ has also
slow recurrence to $\cS$ with respect to $\mu$.  Hence by
Theorem~\ref{thm:tempos-hip-existem} there are
$\sigma,\delta,b>0$ such that $\mu$-almost all $x\in M$
admits infinitely many $(\sigma,\delta,b)$-hyperbolic times
with positive frequency at infinity. Thus there exists a
finite partition $\cP_0$ of $M$ which is generating with
respect to $\mu$.

Indeed let $\cE=\{B(x_i,\delta_1/8), i=1,\dots,l\}$ be a
finite open cover of $M$ by $\delta_1/8$-balls whose
boundary has zero $\mu$ measure. From this we define a
finite partition $\cP_0$ of $M$ as follows.  Start by
setting $P_1=B(x_1,\delta_1/8)$ as the first element of the
partition. Then, assuming that $P_1,\dots, P_k$ are already
defined, set
$P_{k+1}=B(x_{k+1},\delta_1/8)\setminus(P_1\cup\dots\cup
P_k)$ for $k=1,\dots, l-1$. Note that if $P_k\neq\emptyset$
then $P_k$ has non-empty interior, diameter smaller than
$\delta_1/4$ and the boundary $\partial P_k$ is a (finite)
union of pieces of boundaries of balls in a Riemannian
manifold. Thus $\partial P_k$ has zero Lebesgue measure and
zero $\mu$-measure also. Define $\cP_0$ by the elements
$P_k$ constructed above which are non-empty.  Note that
$\mu(\partial\cP_0)=\Leb(\partial\cP_0)=0$ and by the
existence of infinitely many $(\sigma,\delta,b)$-hyperbolic
times for $\mu$-almost every $x$ it is not difficult to see
that $\diam\big(\bigvee_{j=0}^{n-1}f^{-j}\cP_0(x)\big)
\xrightarrow[n\to+\infty]{}0$.

Therefore, since $\mu$ satisfies the Entropy Formula, we can
write
\begin{align*}
  \frac1n\int S_n \log J
  \,d\mu=\mu(\log J)=h_\mu(f)=h_\mu(f,\cP_0)
  \le
  \frac1n H_\mu(\cP_n)
  = \frac1n\int -\log \mu\big(\cP_n(x)\big)\, d\mu
\end{align*}
where $\cP_n=\bigvee_{j=0}^{n-1}f^{-j}\cP_0$ for
$n\ge1$. Hence 
by Jensen's Inequality we get, denoting
$J_n(x)=\prod_{j=0}^{n-1}J\big(f^j(x)\big)$
\begin{align*}
  0\ge\int\log\big[
  J_n(x)\cdot\mu\big(\cP_n(x)\big)\big]\,d\mu(x)
  \ge\log\int J_n(x)\cdot\mu\big(\cP_n(x)\big) \,d\mu(x).
\end{align*}
If we define $Q_\gamma^n=\{x\in M:S_n
J(x)\cdot\mu\big(\cP_n(x)\big)>\gamma\}$ we obtain
\begin{align}\label{eq:Qgamma}
\mu(Q_\gamma^n)\le\gamma^{-1}
\quad\text{for all}\quad
n\ge1.   
\end{align}
Now choose $\gamma_n>0$ such that
$\sum_n\gamma_n^{-1}<\infty$.  Then for $\mu$-almost every
$x\in M$ there exists $n_0\in\NN$ such that for all $n\ge
n_0$ we have $x\not\in Q^n_{\gamma_n}$, i.e.
$J_n(x)\cdot\mu\big(\cP_n(x)\big)\le\gamma_n$ for all
$n\ge n_0=n_0(x)$.
Observe that by the definition and properties of hyperbolic
times, we have that there exists $C_1>0$ such that
\begin{align*}
C_1^{-1}\cdot\Leb\big(\cP_0(f^n(x))\big)
\le
\Leb\big(\cP_n(x)\big)\cdot J_n(x)
\le
C_1 \cdot\Leb\big(\cP_0(f^n(x))\big)  
\end{align*}
whenever $n$ is a hyperbolic time for $x$. This shows that
the $\mu$-measure of the atoms of $\cP_n$ can be bounded
from above by the volume of the same atoms at big enough
hyperbolic times
\begin{align}\label{eq:muLeb}
  \mu\big(\cP_n(x)\big)\le C_0\gamma_n\Leb(\cP_n(x)),
\end{align}
where $C_0=C_1\sup_{x\in M}\Leb\big(\cP_0(x)\big)$. The
hyperbolic times satisfying this condition will be called
\emph{$\mu$-hyperbolic times}.
To use this we need some way to cover any set using
atoms of the sequence $(\cP_n)_n$ at $\mu$-hyperbolic times.

\subsection{Coverings by hyperbolic times}
\label{sec:coverings-hyperb-tim}

Let $\mu$, $f$ and $(\cP_n)_{n\ge0}$ be as in the
previous subsection.  Note that since $f$ is regular and
$\mu$ is $f$-invariant the boundary of $g(P)$ still has zero
Lebesgue measure and zero $\mu$-measure for every atom
$P\in\cP_0$ and every inverse branch $g$ of $f^n$, for any
$n\ge1$.


We can now state the following flexible covering lemma with
$\mu$-hyperbolic preballs. It will enable us to
approximate the $\mu$-measure of a given set through the
measure of families of $\mu$-hyperbolic preballs.

\begin{lemma}[The Hyperbolic Covering Lemma]
  \label{le:coverhiptimes}
  Let a measurable set $E\subset M$, $m\ge1$ and $\zeta>0$
  be given with $\mu(E)>0$. Let $\theta>0$
  be a lower bound for the density of $\mu$-hyperbolic times
  for $\mu$-almost every point. Then there are integers
  $m<n_1<\dots<n_k$ for $k=k(\zeta)\ge1$ and families
  $\cE_i$ of subsets of $M$, $i=1,\dots, k$ such that
\begin{enumerate}
\item $\cE_1\cup\dots\cup\cE_k$ is a finite
  pairwise disjoint family of subsets of $M$;
\item $n_i$ is a $(\sigma/2,\delta/2)$-$\mu$-hyperbolic time
  for every point in $P$, for every element $P\in\cE_i$,
  $i=1,\dots,k$;
\item every $P\in\cE_i$ is the preimage of some element
  $Q\in\cP$ under an inverse branch of $f^{n_i}$, $i=1,\dots,k$;
\item there is an open set $U_1\supset E$ containing the
  elements of $\cE_1\cup\dots\cup\cE_k$ with
  $\mu(U_1\setminus E)<\zeta$;
\item $\mu\Big( E\triangle \bigcup_i \cE_i \Big) \le
  \left(1-\frac{\theta}4 \right)^k < \zeta$.
\end{enumerate}

\end{lemma}

The proof is completely presented in~\cite[Lemma
3.5]{araujo-pacifico2006} and follows~\cite[Lemma
8.2]{OliVi2005} closely.

\subsection{Absolute continuity}
\label{sec:absolute-contin}

We are now ready to deduce that any measure $\mu$ as in the
statement of Theorem \ref{thm:reciprocal-EF} is absolutely
continuous. Indeed observe that, by \eqref{eq:Qgamma} and
the choice of $(\gamma_n)_{n\ge1}$, for any given $\eta>0$
we can find $N=N(\eta)\in\NN$ such that
$\Gamma_\eta=\cap_{n\ge N} \big(M\setminus
Q^n_{\gamma_n}\big)$ satisfies $\mu(\Gamma_\eta)\ge1-\eta$.

Let $E\subset M$ be given with $\mu(E\cap\Gamma_\eta)>0$.
Let $m=N$ in the statement of the Covering
Lemma~\ref{le:coverhiptimes} and set $\zeta>0$ small. Then
we get $\mu\big((E\cap\Gamma_\eta)\triangle \cup_{i=1}^k
\cE_i\big)<\zeta$ where all elements of $\cE_i$ are
$\mu$-preballs and atoms of $\cP_{n_i}$ satisfying the bound
\eqref{eq:muLeb}. In particular by the choice of $m$ we have
$\cup_i\cE_i\subset \Gamma_\eta$ and so we may write
\begin{align}
  \label{eq:abscont}
  \mu(E)=\mu\big(E\cap M\setminus\Gamma_\eta\big)
  + \mu\big(E\cap\Gamma_\eta\big)
  \le
  \eta+\zeta+\mu\big(E\cap\cup_i\cE_i\big)
  \le
  \eta+\zeta+C_0\gamma_{n_k}\Leb\big(E\cap\cup_i\cE_i\big),
\end{align}
where $n_k$ is the largest $\mu$-hyperbolic time used in the
cover given by the Hyperbolic Covering Lemma. 

Hence if we start with a subset $E$ with $\Leb(E)=0$ and
assume that $\mu(E)>0$, then 
there exists $\eta_0$ such that $\mu(E\cap\Gamma_\eta)>0$
for all $0<\eta\le\eta_0$.  Therefore given $\zeta>0$ as
above we obtain \eqref{eq:abscont}. But since $\Leb(E)=0$ we
get $\mu(Z)\le\eta+\zeta$, for all $0<\eta\le\eta_0$, that
is $\mu(Z)\le\zeta$. This is a contradiction since we may
take $\zeta>0$ as small as we like.

We have shown that if $\Leb(E)=0$ then $\mu(E)=0$,
i.e. $\mu\ll\Leb$. Then since the basins of the physical
measures of $f$ cover $M$ except for a volume zero subset,
then it follows easily by the Ergodic Theorem that
$\mu=\sum_{i=1}^k \mu\big(B(\mu_i)\big)\cdot\mu_i$, that is
$\mu\in\overline\co(\FF)$.

Reciprocally, let us now assume that $\mu$ is an
$f$-invariant absolutely continuous probability
measure. Then as above we have $\mu\in\overline\co(\FF)$ and
and thus for some constants $\alpha_i\ge0$ such that
$\sum_i\alpha_i=1$ we have $ h_\mu(f)=\sum_{i=1}^k \alpha_i
h_{\mu_i}(f) = \sum_{i=1}^k \alpha_i\mu_i(\log J) = \mu(\log
J).$ This concludes the proof of
Theorem~\ref{thm:reciprocal-EF}.


\section{Exponentially slow approximation to singularities}
\label{sec:exponent-slow-approx}

Here we apply the (by now standard) arguments of Benedicks
and Carleson, first presented in \cite{BC85,BC91}, to show
that Lorenz-like maps have exponentially slow recurrence to
singularities. This completes the presentation of the
examples in Section~\ref{sec:semifl-over-non-1}.

Let $f:M\to M$ be a one-dimensional $C^{1+\alpha}$ piecewise
expanding map with at most countably many smoothness domains
for some $\alpha\in(0,1)$ as in
Subsection~\ref{sec:semifl-over-non-1}, that is
$|f'|\ge\sigma>1$ and the non-degenerate singular set $\cS$
equals the boundaries of the smoothness domains and
satisfies all the conditions (S1) through (S5).  Then
$\cS=\{b_n\}_n$ where we may assume that the sequence is
strictly monotonous (in counter-clockwise order if
$M=\sS^1$).

We consider the middle points $c_n=(b_n+b_{n+1})/2$ for all
applicable indexes $n$ to define a Lebesgue modulo zero
partition $\cP_0$ of $M$ as follows.  

\subsection{Initial partition}
\label{sec:initial-partit}

Partition $(b_n,c_n)$ into subintervals
\begin{align}
  \label{eq:infty-partition}
M(2n,p)=\big(b_n+d_{2n}e^{-p}, 
b_n+d_{2n}e^{-(p-1)}\big),
\end{align}
where $d_{2n}=c_n-b_n$ and partition the interval
$(c_{n-1},b_n)$ into the following subintervals
\begin{align} \label{eq:infty-partition2}
  M(2n-1,p)=\big(b_n-d_{2n-1}e^{-(p-1)},
  b_n-d_{2n-1}e^{-p}\big)
\end{align}
where $d_{2n-1}=b_n-c_{n-1}$, for all $p\ge1$.  The sets
defined above form a partition of $M$ Lebesgue modulo zero
consisting of small intervals whose length is exponentially
small with respect to the distance to $\cS$.  Let
$\cS'=\cS\cap f(M)$ be the set of singular points of $f$
which matter for the asymptotic dynamics of $f$. 

To define the initial partition consider a threshold
$\rho_0\in\NN$ such that
\begin{align}\label{eq:rho-0}
  e^{-\beta\rho_0}<1 \quad\text{and}\quad
  \left(1+\frac{2}{\rho_0}\right)\left(1+\frac{\rho_0}2\right)^{2/\rho_0}
  < e^\beta
\end{align}
and let $\cP_0$ be formed
by the collection of all intervals $M(n,p)$ for all $n$ and
every $p\ge\rho_0$ together with the connected components of
$M\setminus\big(\cup_{n; p\ge\rho_0}
M(n,p)\cup\{c_n\}_n\big)$, which will be denoted by
$M(n,\rho_0-1)$ whenever they are adjacent to $M(n,\rho_0)$.

For each element $\eta$ of $\cP_0$ denote by $\eta^+$ the
interval obtained by joining $\eta$ with its two neighboring
intervals in $\cP_0$. From \eqref{eq:infty-partition} and
\eqref{eq:infty-partition2} we have the following relations
for all $k$ and every $p\ge \rho_0-1$
\begin{align}\label{eq:+9}
  \Leb(M(k,p)^+)\le 9 \Leb(M(k,p)) = 9
  d_k\cdot e^{-p}(e-1).
\end{align}

\subsection{Refining the partition}
\label{sec:refining-partit}

The partition $\cP_0$ is dynamically refined so that any
pair $x,y$ of points in the same atom of the $n$th
refinement $\cP_n$ belong to the same element of $\cP_0$
during the first consecutive $n$ iterates, i.e.
$\cP_0(f^i(x))=\cP_0(f^i(y))$ for $i=0,\dots,n-1$.  Moreover
$f^n\mid\omega$ is a diffeomorphism for every interval
$\omega\in\cP_n$.

The refinement is defined inductively.  Assume that $\cP_n$
is already defined and for each $\omega\in\cP_n$ there are
sets $R_n(\omega)$ of splitting times and $D_n(\omega)$ of
corresponding splitting depths, to be defined below.

If $f^{n+1}(\omega)$ intersects three or
fewer elements of $\cP_0$, then we set $\omega\in\cP_{n+1}$,
$R_{n+1}(\omega)=R_n(\omega)$ and
$D_{n+1}(\omega)=D_n(\omega)$.  Otherwise consider the
subsets $\eta'=\big(f^{n+1}\mid \omega\big)^{-1}(\eta)$ of the
interval $\omega$, for all elements $\eta$ of $\cP_0$ which
intersect $f^{n+1}(\omega)$.

The family $\{\eta'\}$ obtained above is a partition of
$\omega$. Observe that $f^{n+1}(\eta')$ is either equal to
some $\eta\in\cP_0$ or strictly contained in some
$\eta\in\cP_0$. In the latter case we redefine the partition
joining some of the extreme intervals of $\{\eta'\}$ with
its neighbors so that the new partition $\{\zeta\}$ of
$\omega$ satisfies: for each $\zeta$ there exists
$\eta=M(k,p)\in\cP_0$ such that $\eta\subseteq f^n(\zeta)
\subseteq \eta^+$.

Finally we set $\zeta\in\cP_{n+1}$,
$R_{n+1}(\zeta)=R_n(\zeta)\cup\{n+1\}$ and
$D_{n+1}(\zeta)=D_n(\zeta)\cup\{(k,p)\}$, for each element
of the partition $\{\zeta\}$ of $\omega$ constructed above.
For these elements of $\cP_{n+1}$ we say that $n+1$ is a
\emph{splitting time} and the pairs $(k,p)$ are the
corresponding \emph{splitting depths}.  Repeat the procedure
for each $\omega\in\cP_n$. This completes the construction
of $\cP_{n+1}$ from $\cP_n$ for all $n\ge0$.

\subsection{Bounded distortion}
\label{sec:bounded-distort}

The uniform expansion of length during $n$ iterates ensures
that we have bounded distortion of lengths on atoms of the
partition $\cP_n$.

Indeed let $\omega\in\cP_n$ for some $n\ge1$ and let
$x,y\in\omega$. Note that $f^i\mid\omega$ is a
diffeomorphism for $i=1,\dots,n$, $f$ expands distances at a
minimum rate of $\sigma$ and $f'$ is $\alpha$-H\"older. Then
there exist constants $C,D>0$ such that
\begin{align}
  \log\left|\frac{(f^n)'(x)}{(f^n)'(y)}\right|
  &=
  \sum_{i=0}^{n-1}\big|\log|f'(f^i(x))|-\log|f'(f^j(y)|\big|
  \le
  \sum_{i=0}^{n-1} C \cdot \frac{\big| f^i(x)-f^j(y)
    \big|^\alpha}{\max\{|f'(f^i(x))|, |f'(f^i(y))|\}}
  \nonumber
  \\
  &\le
  \frac{C}\sigma \sum_{i=0}^{n-1}
  \sigma^{i-n}\cdot\big|f^n(x)-f^n(y)\big|^\alpha
  \le D,\label{eq:distortion}
\end{align}
where $D$ depends only on $\sigma$ and on the diameter of $M$.

\subsection{Measure of atoms of $\cP_n$ and return depths}
\label{sec:dist-depth}

Here we show that we can estimate the measure of an element
of $\cP_n$ using the information stored in $R_n$ and $D_n$.

For any given $n\ge1$ and $\omega\in\cP_n$ we have
\begin{itemize}
\item a sequence of times $R_n(\omega)=\{r_1<\dots<r_s\}$
  with $r_1\ge1$ and $r_s\le n$, and
\item a sequence of intervals
  $\omega_0\supsetneq\omega_1\supsetneq\cdots\supsetneq\omega_s=\omega$
  with corresponding depths
  $D_n(\omega)=\{(k_1,p_1),\dots,(k_s,p_s)\}$, where
  $\omega_0\in\cP_0$ and
  $\omega_i\in\cP_{r_i}\cap\dots\cap\cP_{r_{i+1}-1}$
\end{itemize}
such that
\begin{align}\label{eq:encaixe}
  M(k_i,p_i)\subseteq f^{r_i}(\omega_{i-1}) \subseteq M(k_i,p_i)^+
\end{align}
for all $i=0,1,\dots,s-1$ with $r_0=0$ and $s_0=0$. These
times are the iterates where the images of the previous
element of the partition was broken into smaller intervals
as in Subsection~\ref{sec:refining-partit}. Using the
bounded distortion given by~\eqref{eq:distortion} we get
\begin{align*}
  \frac{\Leb(\omega)}{\Leb(\omega_0)} &=
  \frac{\Leb(w_s)}{\Leb(\omega_{s-1})}\cdots
  \frac{\Leb(\omega_1)}{\Leb(\omega_0)}
  \le 
  \prod_{i=1}^s
  D\frac{\Leb\big(f^{r_i}(\omega_{i})\big)}
  {\Leb\big(f^{r_i}(\omega_{i-1})\big)}.
\end{align*}
Now using \eqref{eq:encaixe} and (S5) we bound the last
expression from above by
\begin{align*}
 \prod_{i=1}^s
  \frac{D \Leb\big(M(k_i,p_i)^+\big)}{B^{-1}e^{\beta p_{i-1}}
    d_{k_{i-1}}^{-\beta}(e-1)^{-\beta} \Leb\big(M(k_{i-1},p_{i-1})\big)}
\end{align*}
and using \eqref{eq:+9} this can be easily simplified
yielding
\begin{align}\label{eq:omega0}
  \Leb(\omega)
  \le
  \prod_{i=0}^{s-1} d_{k_i}^\beta  e^{-2\beta p_i}
  \le
  \exp\left(-\beta\sum_{i=0}^{s-1}(p_i+q_i)\right)
\end{align}
where $q_i=[-\log d_{k_i}]$ with $[z]=\max\{k\in\ZZ:k\le
z\}$. We have used $p\ge\rho_0$ and
$\log(9BD(e-1)^\beta)/\rho_0\le\beta$ to compensate the
constants on the exponent $2\beta$. Recall also that
$\omega_0=M(k_0,p_0)$. Note also that if
$R_n(\omega)=\emptyset$, then since there is no splitting but
there is uniform expansion together with distortion control,
we get
\begin{align}\label{eq:omega1}
  \Leb\big(f^n(\omega)\big)=\int_{\omega}\big|(f^n)\big|\,d\Leb
  \ge D\sigma^n\Leb(\omega)
  \text{ so }
  \Leb(\omega)\le D^{-1}\sigma^{-n}.
\end{align}

\subsection{Distance to $\cS$ and splitting depths}
\label{sec:distance-cs-splitt}

Let again $n\ge1$ and $\omega\in\cP_n$ be given and consider
the sets $R_n(\omega)$ and $D_n(\omega)$.  Consider the
intervals
$\omega_0\supsetneq\omega_1\supsetneq\cdots\supsetneq\omega_s=\omega$
as before.  Note that for the iterates $i$ between two
consecutive times $r<r'$ from $R_n$, i.e. if $r<i<r'$ then
there exists $M(l_i,q_i)\in\cP_0$ such that
$f^i(\omega_r)\subseteq M(l_i,q_i)^+$ by this choice of $i$.
Moreover by condition (S5) and by \eqref{eq:infty-partition}
 and \eqref{eq:infty-partition2} we deduce
\begin{align*}
  9d_{l_1}(e-1)e^{-q_1}
  &\ge\Leb\big(f^{r+1}(\omega_r)\big)
  \ge
  \big(B d_{k_r}e^{-p_r}\big)^{-\beta}\Leb\big(f^r(\omega_r)\big)
  \\
  &= 
  \left(
    \frac{B}{e-1}\Leb\big(f^r(\omega_r)\big)\right)^{-\beta}
  \Leb\big(f^r(\omega_r)\big)
  \\
  &=
  \left( \frac{e-1}{B} \right)^\beta 
  \Leb\big(f^r(\omega_r)\big)^{1-\beta}.
\end{align*}
Hence $d_{l_1}e^{-q_1-1}\ge \big(9e(e-1)\big)^{-1}
\Leb\big(f^{r+1}(\omega_r)\big)$ is the estimate for the
minimum distance from $\cS$ to $f^{r+1}(\omega_r)$.  Let
$L_i=\Leb\big(f^{r+i}(\omega_r)\big)$ and
$D_i=\dist\big(f^{r+i}(\omega_r),\cS\big)$ for
$i=0,\dots,r'-r-1$. Then the reasoning above shows that $
L_1\ge \left( \frac{e-1}{B} \right)^\beta L_0^{1-\beta}$ and
$L_{i+1}\ge \left( \frac{9e(e-1)}{B} \right)^\beta
L_i^{1-\beta}$, and also $D_i\ge L_i/\big(9e(e-1)\big)$ for
$i=1,\dots,r'-r-1$.  It is now easy to see that
\begin{align*}
  -\log L_{i+1}
  &\le
  -(1-\beta)\log L_i+\beta\log\frac{9e(e-1)}{B}
  \\
  &=
  -\left(1-\beta- \beta \log\Big(\frac{9e(e-1)}{B}\Big)
    \Big/ 
    \log L_i\right)\log L_i
  =-\gamma\log L_i
\end{align*}
where we may assume that $\gamma\in(0,1)$ since it is no
restriction to increase the value of $B$ if needed.
Hence
\begin{align*}
  -\!\!\!\sum_{i=1}^{r'-r-1} \log\dist\big(f^{r+i}(\omega_r),\cS\big)
  &\le
  -\!\!\!\sum_{i=1}^{r'-r-1}\Big(\log
  L_i-\log\big(9e(e-1)\big)\Big)
  \\
  &\le
  -\mathrm{Const}\cdot\log L_0 + (r'-r)
  \log\big(9e(e-1)\big)
  \\
  &\le 
  -\mathrm{Const}\cdot\log L_0,
\end{align*}
since by uniform expansion and by definition of $r'$ we have
$\sigma^{r'-r}L_0\le1$ and also
$r'-r\le-\log(L_0)/\log\sigma$. Since $r<r'$ were two
arbitrary consecutive elements of $R_n(\omega)$ for
$\omega\in\cP_n$ we have shown that
\begin{align}  \label{eq:aproxdepth}
  \sum_{j=0}^{s-1}-\log\dist\big(f^j(x),\cS\big)\le
  -\mathrm{Const}\sum_{(k,p)\in D_s(\omega)}\log\big(d_ke^{-p}\big)
\end{align}
for all $x\in\omega$, where $s<n$ is the last splitting time
before $n$ ($s=\max R_n(\omega)$). 

However if $m>n$ is the first integer such that
$\omega\not\in\cP_m$ but $\omega\in\cP_{l}$ for $n<l<m$,
then we can write the following disjoint union
$\omega=\bigcup_{\omega'\in\cP_m} \omega'\cap\omega$.
Repeating the argument for $x\in\omega'\cap\omega$ for each
$\omega'\in\cP_m$ intersecting $\omega$ we can obtain a
relation like~\eqref{eq:aproxdepth} with $D_s(\omega)$
replaced by $D_n(\omega)$ as the summation range, where $n$
is between $s$ and $m$. This shows that \emph{the average of
  the $\log$ of the distance to the singular set is bounded
  by the sum of the depths at splitting times modulo a
  constant}.

\subsection{Expected value of splitting depths}
\label{sec:expect-value-splitt}

Now we estimate the expected value of the splitting depths
for deep splitting times up to $n$ iterates of the map.
Define for a co-countable set
of $x\in M$ the function $\D_n(x)=-\sum_{(k,p)\in
  D_n(\cP_n(x))}\log(d_ke^{-p})$ where $\cP_n(x)$ is the
unique atom of $\cP_n$ which contains $x\in M$. Define also
the truncated sum: for any given $\delta>0$ set for the same
points $x\in M$ as above
\begin{align}\label{eq:D_n-delta}
\D_n^\delta(x)=\sum_{\substack{(k,p)\in
  D_n(\cP_n(x))\\ d_k e^{-p}<\delta}}  -\log(d_ke^{-p}).
\end{align}
By the arguments in Subsection~\ref{sec:distance-cs-splitt}
and by the definitions \eqref{eq:infty-partition} and
\eqref{eq:infty-partition2} we obtain
\begin{align}
  \label{eq:logdist-Dn}
\sum_{j=0}^{n-1}-\log\dist_\delta\big(f^j(x),\cS\big)
  \le \D_n^\delta(x).
\end{align}
Define the \emph{number of splittings up to the $n$th iterate}
$t_n(x)=\# R_n\big(\cP_n(\omega)\big)$ and also the \emph{number of
deep splittings among these}
$u_n(x)=\#\big\{(k,p)\in
R_n\big(\cP_n(\omega)\big): d_k e^{-p}<\delta \big\}.$

Given $x$ and $n\ge1$ we let $0=r_0<r_1<\dots<r_t$ with
$t=t_n(x)$ be the splitting times along the orbit of $x$ up
to the $n$th iterate and $0\le s_1<\dots<s_u$ be indexes
corresponding to deep splittings, where $u=u_n(x)$ in what
follows. Note that each quantity above is constant on the
elements of $\cP_n$. Define
\begin{align*}
  A^{u,t}_{(\kappa_1,\rho_1),\dots,(\kappa_u,\rho_u)}(n)
  =
  \big\{
  x\in M: t_n(x)=t, \, u_n(x)=n\,\text{ and }
  (k_{s_i},p_{s_i})=(\kappa_i,\rho_i),\, i=1,\dots,u
  \big\}
\end{align*}
the set of points which in $n$ iterates have $t$ splitting
times and $u$ deep splittings among these, with the specified
depths $(\kappa_1,\rho_1),\dots,(\kappa_u,\rho_u)$.

\begin{lemma}
  \label{le:A-u-t}
  $\Leb\big(
  A^{u,t}_{(\kappa_1,\rho_1),\dots,(\kappa_u,\rho_u)}(n)\big)
  \le \binom{t}{u}
  \exp\left(-\beta\sum_{i=1}^{u}(\eta_i+\rho_i)\right)$
  where $\eta_i=[-\log d_{\kappa_i}]$.
\end{lemma}

\begin{proof}
  Using the estimate \eqref{eq:omega0} we get the following
  bound for the Lebesgue measure of $
  A^{u,t}_{(\kappa_1,\rho_1),\dots,(\kappa_u,\rho_u)}(n)$
  \begin{align*}
  \binom{t}{u}
  &\exp\left(-\beta\sum_{i=1}^{u}(\eta_i+\rho_i)\right)
  \cdot
  \exp\Big(-\beta\!\!\!\!\!\!
    \sum_{\substack{(k_j,p_j)\text{
          s.t. }d_ke^{-p_j}\ge\delta
        \\j=1,\dots,t-u}}
    \!\!\!\!\!\!(\nu_j+p_j)\Big).
  \end{align*}
  The binomial coefficient takes into account all the
  possible orderings of sequences of $u$ deep splitting
  times among $t$ splitting times and the last exponential
  bounds the contribution of all the possible $t-u$ non-deep
  splitting times, with $\nu_j=[-\log d_{k_j}]$. But since
  $p\ge\rho_0$ was chosen as in \eqref{eq:rho-0} and
  $\nu_j\ge0$ we conclude that the last exponential is
  smaller than 1. So we obtain the bound in the statement.
\end{proof}

\begin{lemma}
  \label{le:esperado}
  For any $z>\beta$ we have
  $\int e^{z \D_n^\delta(x)} \, dx
    \le
    e^{\theta(\delta)n}$
  where $\theta(\delta)$ is such that $\theta(\delta)\searrow0$ when
  $\delta\searrow0$.
\end{lemma}

\begin{proof}
By definition 
\begin{align}
  \int e^{z \D_n^\delta(x)} \, dx
  &=
  \sum_{\omega\in\cP_n}e^{z\D_n^\delta(\omega)}\cdot\Leb(\omega)
  \le \sum_{\substack{\omega_0\in\cP_0\\D\sigma^n\Leb(\omega_0)\le1}}
    \Leb(\omega_0)\nonumber
  \\
  &\quad
  +\sum_{0< u\le t<n}
  \sum_{\substack{(\kappa_i,\rho_i)\\i=1,\dots,u}} 
  e^{z\D_n^\delta(\omega)}
  \Leb\big( A^{u,t}_{(\kappa_1,\rho_1),\dots,(\kappa_u,\rho_u)}(n)\big)
  \label{eq:somona}
\end{align}
where we are considering all possible combinations of
splitting depths and of deep splittings among all the
splitting times, for all elements of $\cP_n$ in the second sum.

Consider the first term corresponding to the atoms of
$\cP_0$ which were not split during the first $n$ iterates.
This sum can be separated as follows
\begin{align}
  \sum_{\substack{\omega_0\in\cP_0\\D\sigma^n\Leb(\omega_0)\le1}}
  \hspace{-0.6cm}\Leb(\omega_0)
  =
  \sum_{\substack{D\sigma^n\Leb(\omega_0)\le1
    \\
    d_k\sigma^{n/2}<1}} \hspace{-0.6cm}\Leb(\omega_0)
  +
  \hspace{-0.6cm}\sum_{\substack{D\sigma^n\Leb(\omega_0)\le1
    \\
    d_k\sigma^{n/2}\ge1}} \hspace{-0.6cm}\Leb(\omega_0)
  \le
  \Leb\big( B(\cS,\sigma^{-n/2})\big)
  +\hspace{-1cm}
  \sum_{p>\log\big(D(e-1)\sigma^{n/2}\big)}\hspace{-1cm}e^{-p}
  \le C e^{-c n}\label{eq:soma1}
\end{align}
for some constants $C,c>0$, where we have used
expression~\eqref{eq:+9} for the length of the atoms of
$\cP_0$ in terms of $(k,p)$ together with condition (S4) and
the obvious $d_k>0$ and $\sum_k d_k=1$. Note that if $\cS$
is finite then the condition $d_k\sigma^{n/2}<1$ is always
false for big enough $n$. So in this case we only have the
right hand side sum above.

Now we bound the second term~\eqref{eq:somona}. 
Considering Lemma~\ref{le:A-u-t} and taking into account
$\D_n^\delta$ we obtain (with $\eta_j=[-\log d_{\kappa_j}]$)
\begin{align*}
  \sum_{0< u\le t<n}
  \sum_{\substack{(\kappa_i,\rho_i)\\i=1,\dots,u}} \!\!\!
  \binom{t}{u} e^{-(\beta+z)\sum_i (\eta_i+\rho_i)}
  \le\!\!\!
  \sum_{0< u\le t<n}\sum_{h>u\ell(\delta)}
  \!\!\!\binom{t}{u} 
  u L(h,u) e^{-(\beta+z)h}
\end{align*}
where $h=\sum_i(\eta_i+\rho_i)$, $\ell(\delta)$ is an
integer such that every pair $(k,p)$ satisfying
$d_ke^{-p}<\delta$ also satisfies $k+p>\ell(\delta)$,
and 
\begin{align*}
  L(h,u)=\#\left\{
    \big( (\eta_i,\rho_i)\big)_{i=1,\dots,u}\in\NN_0^{2u}:
    \sum_{i=1}^{u}(\eta_i+\rho_i)=h \text{ with }
    \rho_i\ge\rho_0
  \right\}.
\end{align*}
Moreover the factor $u$ bounds the number of distinct
$d_{k_i}$ with the same value $\eta_i$ along the $n$
iterates of the orbit of the points. Observe that
\begin{align*}
  L(h,u)
  \le\#\left\{
    (h_i)\in\NN_0^{2u}: \sum_{i=0}^{2u}h_i=h
    \right\}
    =
    \binom{h+2u-1}{2u-1}
\end{align*}
and by a standard application of Stirling's Formula
\begin{align*}
  L(h,n)
  \le
  \left(
    c^{1/h}\big(1+\frac{2u-1}{h}\big)
    \big(1+\frac{h}{2u-1}\big)^{(2u-1)/h}
    \right)^h
    \le e^{\beta h}\le e^{zh}
\end{align*}
where $0<c<1$ is a constant independent of the other
variables and the last inequalities follow by $h\ge \rho_0 u$,
 by the choice of $\rho_0$ in \eqref{eq:rho-0} and by taking
 $z>\beta$.

Collecting the bounds we have obtained we conclude that the
second sum in \eqref{eq:somona} can be bounded by the
following expression
\begin{align*}
  \sum_{0< u\le t<n} \binom{t}{u} 
  u \sum_{h>u\ell(\delta)} e^{-\beta h}
  &\le
  \sum_{u=0}^n n \binom{n-1}{u} \cdot
  u e^{-\beta u \ell(\delta)/2} \cdot \frac{e^{-\beta u
    \ell(\delta)/2}}{1-e^{-\beta}}
  \le
  \sum_{u=0}^n \binom{n}{u} C 
  \frac{\big(e^{-\beta\ell(\delta)/2}\big)^u}{1-e^{-\beta}}
  \\
  &\le
  \Big( 1 + \frac{C}{1-e^{-\beta}}e^{-\beta\ell(\delta)/2}\Big)^n
\end{align*}
for some constant $C>0$ bounding $\{ u e^{-\beta u
  \ell(\delta)/2}\}_{u\ge0}$ (which can be taken
independently of $\ell(\delta)$).  Finally since
$\ell(\delta)$ grows without limit when $\delta\searrow0$,
the statement of the lemma follows just by increasing the
value of $C$ to take into account the small bound of the
first sum~\eqref{eq:soma1}.
\end{proof}


\subsection{Measure of the points with bad recurrence}
\label{sec:measure-points-with}

We are now ready to deduce exponentially slow approximation
to the singular set $\cS$. Indeed we just have to use
Tchebishev's inequality, as follows: given
$\epsilon,\delta>0$ we know there exists a constant $C>0$ as
in Subsection~\ref{sec:distance-cs-splitt} such that
\begin{align*}
  \Big\{x\in
    M:-\frac1n\sum_{i=0}^{n-1}\log\dist_\delta\big(f^i(x),\cS\big)
    \ge\epsilon\Big\}
  \subseteq
  \Big\{x:
  \frac{\D_n^\delta(x)}{n}\ge\frac{\epsilon}{C}\Big\}
  =
  \Big\{x:e^{z\D_n^\delta(x)}\ge e^{n\epsilon/C}\Big\}
\end{align*}
hence
\begin{align*}
  \Leb\Big\{x\in
    M:-\frac1n\sum_{i=0}^{n-1}\log\dist_\delta\big(f^i(x),\cS\big)
    \ge\epsilon\Big\}
  \le
    e^{-n\epsilon/C}\int e^{z\D_n^\delta}\,d\Leb 
  =
    e^{-n\big(\epsilon/C-\theta(\delta)\big)}
\end{align*}
which can be made exponentially small by choosing $\delta>0$
small enough so that $\epsilon/C>\theta(\delta)$. This
proves that a piecewise expanding map $f$ in our settings
has exponentially slow recurrence to the singular set,
completing the proof of the statements in
Subsection~\ref{sec:semifl-over-non-1} and of
Corollary~\ref{mcor:LDLorenz} after the reduction procedure
of Subsection~\ref{sec:singul-hyperb-attrac}.



\def\cprime{$'$}

\end{document}